\documentclass[a4paper]{amsart}
\usepackage{lmodern}

\usepackage[T1]{fontenc}
\usepackage[latin9]{inputenc}
\usepackage[english]{babel}
\usepackage{amsmath}
\usepackage{amssymb}

\usepackage[notcite]{showkeys}
\usepackage{hyperref}
\begin{document}
\noindent \begin{center}
\textbf{\LARGE Saccharinity with ccc}
\par\end{center}{\LARGE \par}

\noindent \begin{center}
{\large Haim Horowitz and Saharon Shelah}
\par\end{center}{\large \par}

\noindent \begin{center}
\textbf{\small Abstract}
\par\end{center}{\small \par}

\noindent \begin{center}
{\small Using creature technology, we construct families of Suslin
ccc non-sweet forcing notions $\mathbb Q$ such that $ZFC$ is equiconsistent
with $ZF+$''Every set of reals equals a Borel set modulo the $(\leq \aleph_1)$-closure
of the null ideal associated with $\mathbb Q$''$+$''There is an
$\omega_1$-sequence of distinct reals''. This answers a question of the second author and Kellner. As an application of independent interest, we also show how our forcing adds a new $\Pi^1_2$ singleton over $L$ without relying on $L$-combinatorics.\footnote{{\small Date: May 26, 2025}{\small \par}



Publication 1067 of the second author%
}}
 \par\end{center} {\small \par}

\textbf{\large 1. Introduction}{\large \par}

\textbf{Some history}

The study of the consistency strength of regularity properties originated
in Solovay's celebrated work {[}So2{]}, where he proved the following
result:

\textbf{Theorem ({[}So2{]}): }Suppose there is an inaccessible cardinal,
then after forcing (by Levy collapse) there is an inner model of $ZF+DC$
where all sets of reals are Lebesgue measurable and have the Baire
property.

Following Solovay's result, it was natural to ask whether the existence
of an inaccessible cardinal is necessary for the above theorem. This
problem was settled by Shelah ({[}Sh176{]}) who proved the following
theorems:

\textbf{Theorem ({[}Sh176{]}): }1. If every $\Sigma_3^1$ set of reals
is Lebesgue measurable, then $\aleph_1$ is inaccessible in $L$.

2. $ZF+DC+"$all sets of reals have the Baire property$"$ is equiconsistent
with $ZFC$.

A central concept in the proof of the second theorem is the amalgamation
of forcing notions, which allows the construction of a suitably homogeneous
forcing notion, thus allowing the use of an argument similar to the
one used by Solovay, in which we have {}``universal amalgamation''
(for years it was a quite well known problem). As the problem was
that the countable chain condition is not necessarily preserved by
amalgamation, Shelah isolated a property known as {}``sweetness'',
which implies $ccc$ and is preserved under amalgamation. See more
on the history of the subject in {[}RoSh672{]}.

\textbf{2. General regularity properties}

Given an ideal $I$ on the reals, we say that a set of reals $X$
is $I-$measurable if $X\Delta B \in I$ for some Borel set $B$,
this is a straightforward generalization of Lebesgue measurability
and the Baire property. 

Given a definable forcing notion $\mathbb Q$ adding a generic real
$\underset{\sim}{\eta}$ (we may write $\mathbb Q$ instead of $(\mathbb Q,\underset{\sim}{\eta})$)
and a cardinal $\aleph_0 \leq \kappa$, there is a natural ideal on
the reals $I_{\mathbb{Q},\kappa}$ associated to $(\mathbb{Q},\kappa)$
(see definition 18), such that, for example, $I_{Cohen, \aleph_0}$
and $I_{Random,\aleph_0}$ are the meagre and null ideals, respectively.
Hence in many cases the study of ideals on the reals corresponds to
the study of definable forcing notions adding a generic real. On the
study of ideals from the point of view of classical descriptive set
theory, see {[}KeSo{]} and {[}So1{]}. For a forcing theoretic point
of view, see {[}RoSh672{]}. Another approach to the subject can be
found in {[}Za{]}.

We are now ready to formulate the first approximation for our general
problem:

\textbf{Problem: }Classify the definable $ccc$ forcing notions according
to the consistency strength of $ZF+DC+"$all sets of reals are $I_{\mathbb{Q},\kappa}-$measurable$"$.

Towards this we may ask: Given a definable $ccc$ forcing notion $\mathbb Q$,
is it possible to get a model where all sets of reals are $I_{\mathbb{Q},\kappa}-$measurable
without using an inaccessible cardinal and for non-sweet forcing notions?

\textbf{3. Saccharinity}

A positive answer to the last question was given by Kellner and Shelah
in {[}KrSh859{]} for a proper \textbf{non-ccc} (very non-homogeneous)
forcing notion $\mathbb{Q}$, where the ideal is $I_{\mathbb{Q},\aleph_1}$.

In this paper we shall prove a similar result for a \textbf{ccc }forcing
notion, omitting the $DC$ but getting an $\omega_1$-sequence of
distinct reals. By {[}Sh176{]}, the existence of such sequence is
inconsistent with the Lebesgue measurability of all sets of reals,
hence our forcing notions are, in a sense, closer to Cohen forcing
than to Random real forcing.

Our construction will involve the creature forcing techniques of {[}RoSh470{]}
and {[}RoSh628{]}, and will result in definable forcing notions $\mathbb{Q}_{\mathbf n}^i$
which are non-homogeneous in a strong sense: Given a finite-length
iteration of the forcing, the only generic reals are those given explicitly
by the union of trunks of the conditions that belong to the generic
set. 

The homogeneity will be achieved by iterating along a very homogeneous
(thus non-wellfounded) linear order. By moving to a model where all
sets of reals are definable from a finite sequence of generic reals,
we shall obtain the consistency of $ZF+"$all sets of reals are $I_{\mathbb{Q}_{\mathbf n}^i,\aleph_1}-$measurable$"+"$There
exists an $\omega_1$-sequence of distinct reals$"$.

It's interesting to note that our model doesn't satisfy $AC_{\aleph_0}$,
thus leading to a finer version of the problem presented earlier:

\textbf{Problem: }Classify the definable $ccc$ forcing notions according
to the consistency strength of $T+"$all sets of reals are $I_{\mathbb{Q},\kappa}-$measurable''
where $T\in \{ZF,ZF+AC_{\aleph_0}, ZF+DC,ZF+DC(\aleph_1), ZFC\}$,
and similarly for $T'=T+WO_{\omega_1}$ where $T$ is as above and
$WO_{\omega_1}$ is the statement $"$There is an $\omega_1$-sequence
of distinct reals$"$.

Remark: Note that for some choices of $T$, $\mathbb Q$ and $\kappa$,
the above statement might be inconsistent.

We intend to address this problem in {[}F1424{]} and other continuations.

\textbf{4. On the special properties of $\mathbb{Q}^2_{\mathbf n}$}

We shall focus in this paper on two types of forcing notions, namely $\mathbb{Q}^1_{\mathbf n}$ and $\mathbb{Q}_{\mathbf n}^2$. In the years following the initial posting of this paper online, the forcing $\mathbb{Q}_{\mathbf n}^2$ was popularized by [KST], where it was reintroduced under the name $\tilde{\mathbb E}$ and shown to play an important role in the study of Cichon's maximum. See [Me1] for a systematic presentation and proofs that the forcing has strong FAM limits and ultrafilter limits for intervals. Another attractive feature of $\mathbb{Q}_{\mathbf n}^2$ which we shall investigate in the end of this paper is the fact that it provides a novel way to add a $\Pi^1_2$ singleton over $L$, without relying on $L$-combinatorics (which was crucial for Jensen's proof).

A remark on notation: 1. Given a tree $T\subseteq \omega^{<\omega}$
and a node $\eta \in T$, we shall denote by $T^{[\eta \leq]}$ the
subtree of $T$ consisting of the nodes $\{ \nu : \nu \leq \eta \vee \eta \leq \nu \}$.

2. For $T$ as above, if $\eta \in T$ is the trunk of $T$, let $T^+:=\{\nu \in T : \eta \leq \nu\}$.

\textbf{\large 2. Norms, $\mathbb{Q}_{\mathbf n}^{1}$ and $\mathbb{Q}_{\mathbf n}^{2}$}{\large \par}

In this section we shall define a collection $\mathbf N$ of parameters.
Each parameter $\mathbf n \in \mathbf N$ consists of a subtree with finite
branching of $\omega^{<\omega}$ with a rapid growth of splitting
and a norm on the set of successors of each node in the tree. 

From each parameter $\mathbf n \in \mathbf N$ we shall define two forcing
notions, $\mathbb{Q}_{\mathbf n}^1$ and $\mathbb{Q}_{\mathbf n}^2$.
We shall prove that they're nicely definable ccc. We will show additional
nice properties in the case of $\mathbb{Q}_{\mathbf n}^2$, such as
a certain compactness property and the fact that being a maximal antichain
is a Borel property. We refer the reader to {[}RoSh470{]} and {[}RoSh628{]}
for more information on creature forcing.

\textbf{Definition} \textbf{1}: 1. A norm on a set $A$ is a function
assigning to each $X\in P(A)\setminus \{\emptyset\}$ a non-negative
real number such that $X_1\subseteq X_2 \rightarrow nor(X_1)\leq nor(X_2)$.

2. Let $\mathbf M$ be the collection of pairs $(\mathbb Q,\underset{\sim}{\eta})$
such that $\mathbb Q$ is a Suslin ccc forcing notion and $\underset{\sim}{\eta}$
is a $\mathbb Q$-name of a real.

\textbf{Definition} \textbf{2}: Let $\mathbf{N}$ be the set of tuples
$\mathbf{n}=(T,nor,\bar{\lambda},\bar{\mu})=(T_{\mathbf n},nor_{\mathbf n},\bar{\lambda}_{\mathbf n},\bar{\mu}_{\mathbf n})$
such that:

a. $T$ is a subtree of $\omega^{<\omega}$.

b. $\bar{\mu}=(\mu_{\eta} : \eta \in T)$ is a sequence of non-negative
real numbers.

c. $\bar{\lambda}=(\lambda_{\eta} : \eta \in T)$ is a sequence of
pairwise distinct non-zero natural numbers such that:

1. $\lambda_{\eta}=\{m : \eta{}^\frown m \in T\}$, so $T\cap \omega^n$
is finite and non-empty for every $n$.

2. If $lg(\eta)=lg(\nu)$ and $\eta<_{lex} \nu$ then $\lambda_{\eta} \ll \lambda_{\nu}$, where $"m \ll k"$ means that "$k$ is much larger than $m$", for our purposes it suffices to require that $m \ll k \iff \beth_m \leq k$.

3. If $lg(\eta)<lg(\nu)$ then $lg(\eta) \ll \lambda_{\eta}\ll \lambda_{\nu}$.

4. $lg(\eta) \ll \mu_{\eta} \ll \lambda_{\eta}$ for $\eta \in T$.

d. For $\eta \in T$, $nor_{\eta}$ is a function with domain $\mathcal{P}^{-}(suc_T(\eta))=\mathcal{P}(suc_T(\eta))\setminus \emptyset$
and range $\subseteq \mathbb{R}_{\geq 0}$ such that: 

1. $nor_{\eta}$ is a norm on $suc_T(\eta)$ (see definition 1).

2. $(lg(\eta)+1)^2 \leq \mu_{\eta}\leq nor_{\eta}(suc_T(\eta))$.

e. $\lambda_{<\eta}:= \Pi \{\lambda_{\nu} : \lambda_{\nu}<\lambda_{\eta}\} \ll \mu_{\eta}$.

f. (Co-Bigness) If $k\in \mathbb{R}^+$, $a_i \subseteq suc_{T_{\mathbf n}}(\eta)$
for $i<i(*)\leq \mu_{\eta}$ and $k+\frac{1}{\mu_{\eta}}\leq nor_{\eta}(a_i)$
for every $i<i(*)$, then $k\leq nor_{\eta}(\underset{i<i(*)}{\cap}a_i)$.

g. If $1\leq nor_{\eta}(a)$ then $\frac{1}{2}<\frac{|a|}{|suc_{T_{\mathbf n}}(\eta)|}$.

h. If $k+ \frac{1}{\mu_{\eta}}\leq nor_{\eta}(a)$ and $\rho \in a$, then $k\leq nor_{\eta}(a\setminus \{\rho\})$.

\textbf{Definition} \textbf{3}: A. For $\mathbf n \in \mathbf N$ we shall
define the forcing notions $\mathbb{Q}_{\mathbf n}^1 \subseteq \mathbb{Q}_{\mathbf n}^{\frac{1}{2}} \subseteq \mathbb{Q}_{\mathbf n}^0$
as follows:

1. $p\in \mathbb{Q}_{\mathbf n}^0$ iff for some $tr(p) \in T_{\mathbf n}$
we have:

a. $p$ or $T_p$ is a subtree of $T_{\mathbf n}^{[tr(p)\leq]}$ (so
it's closed under initial segments) with no maximal node.

b. For $\eta \in lim(T_p)$, $lim(nor_{\eta \restriction l}(suc_{T_p}(\eta \restriction l)) : lg(tr(p)) \leq l<\omega)=\infty$.

c. $2-\frac{1}{\mu_{tr(p)}} \leq nor(p)$ (where $nor(p)$ is defined
in C(b) below).

2. $p\in \mathbb{Q}_{\mathbf n}^{\frac{1}{2}}$ if $p\in \mathbb{Q}_{\mathbf n}^0$
and $nor_{\eta}(Suc_p(\eta))>2$ for every $tr(p) \leq \eta \in T_p$.

We shall prove later that $\mathbb{Q}_{\mathbf n}^{\frac{1}{2}}$ is
dense in $\mathbb{Q}_{\mathbf n}^0$.

3. $p\in \mathbb{Q}_{\mathbf n}^1$ if $p\in \mathbb{Q}_{\mathbf n}^0$
and for every $n<\omega$, there exists $k^p(n)=k(n)>lg(tr(p))$ such
that for every $\eta \in T_p$, if $k(n) \leq lg(\eta)$ then $n\leq nor_{\eta}(Suc_p(\eta))$.

B. $\mathbb{Q}_{\mathbf n}^i \models p\leq q$ $(i\in \{0,\frac{1}{2},1\})$
iff $T_q \subseteq T_p$.

C. a. For $i\in \{0,\frac{1}{2},1\}$, $\underset{\sim}{\eta_{\mathbf n}^{i}}$
is the $\mathbb{Q}_{\mathbf n}^i-$name for $\cup \{tr(p) : p\in \underset{\sim}{G_{\mathbb{Q}_{\mathbf n}^i}}\}$.

b. For $i\in \{0,\frac{1}{2},1\}$ and $p\in \mathbb{Q}$ let $nor(p):=sup \{a\in \mathbb{R}_{>0} : \eta \in T_p^+ \rightarrow a\leq nor_{\eta}(suc_{T_p}({\eta}))\}=inf \{nor_{\eta}(suc_{T_p}(\eta)) : \eta \in T_p\}$.

D. For $i\in \{0,\frac{1}{2},1\}$ let $\mathbf{m}_{\mathbf n}^{i}=\mathbf{m}_{i,\mathbf n}=(\mathbb{Q}_{\mathbf n}^{i}, \underset{\sim}{\eta_{\mathbf n}^{i}})$.

We shall now describe a concrete construction of some $\mathbf n \in \mathbf N$:

\textbf{Definition} \textbf{4}: We say $\mathbf n \in \mathbf N$ is special
when:

a. For each $\eta \in T_{\mathbf n}$ the norm $nor_{\eta}$ is defined
as follows: for $\emptyset \neq a\subseteq suc_T(\eta)$, $nor_{\eta}(a)=\frac{log_*(|suc_T(\eta)|)}{\mu_{\eta}^2}-\frac{log_*|suc_T(\eta)\setminus a|}{\mu_{\eta}^2}$
where $log_*(x)=max \{n: \beth_n \leq x\}$ (where $\beth_0=1$ and $\beth_{n+1}=2^{\beth_n}$).

b. $\mu_{\eta}=nor_{\eta}(suc_{t_{\mathbf n}}(\eta))$.

\textbf{Observation 4A:} There are $T_{\mathbf n}$, $(\lambda_{\eta}, \mu_{\eta} : \eta \in T_{\mathbf n})$
and $(nor_{\eta} : \eta \in T_{\mathbf n})$ satisfying the requirements
of definition 2, where the norm is defined as in definition 4 (hence
$\mathbf n \in \mathbf N$ is special).

\textbf{Proof: }It's easy to check that the following $(T_{\mathbf n}, (\mu_{\eta}, \lambda_{\eta} : \eta \in T_{\mathbf n}))$
together with the norm from definition 4 form a special $\mathbf n \in \mathbf N$
where $T_{\mathbf n} \cap \omega^n$, $(\mu_{\eta}, \lambda_{\eta} : \eta \in T_{\mathbf n} \cap \omega^n)$
are defined by induction on $n<\omega$ as follows: 

a. $T_\mathbf n \cap \omega^0 =\{<>\}$.

b. At stage $n+1$, for $\eta \in T_{\mathbf n} \cap \omega^n$, by
induction according to $<_{lex}$, define $\mu_{\eta}=\beth_{\lambda_{<\eta}}$,
$\lambda_{\eta}=\beth_{{\mu_{\eta}}^3}$ and the set of successors
of $\eta$ in $T_{\mathbf n}$ is defined as $\{\eta{}^\frown (l) : l<\lambda_{\eta}\}$.

For example, we shall prove the co-bigness property:

Suppose that $\eta \in T_{\mathbf n}$ $(a_i : i<i(*))$ are as in definition
2(f). Denote $k_1=|suc_{T_{\mathbf n}}(\eta)|$ and $k_2=max \{|suc_{T_{\mathbf n}}(\eta) \setminus (a_i)| : i<i(*)\}$.
Therefore, $\frac{log_*(k_1)}{\mu_n^2}-\frac{log_*(k_2)}{\mu_n^2}\leq nor_{\eta}(a_i)$
(so necessarily $k+\frac{1}{\mu_{\eta}} \leq \frac{log_*(k_1)}{\mu_n^2}-\frac{log_*(k_2)}{\mu_n^2}$).
Let $a=\underset{i<i(*)}{\cup}a_i$ and $k_3=|suc_{T_{\mathbf n}}(\eta)\setminus a| \leq i(*)k_2 \leq \mu_{\eta}k_2$.
Therefore $\frac{log_*(k_1)}{\mu_{\eta}^2}-\frac{log_*(\mu_{\eta}k_2)}{\mu_n^2} \leq \frac{log_*(k_1)}{\mu_{\eta}^2}-\frac{log_*(k_3)}{\mu_{\eta}^2}=nor_{\eta}(a)$.
We have to show that $k\leq nor_{\eta}(a)$, so it's enough to show
that $k\leq \frac{log_*(k_1)}{\mu_{\eta}^2}-\frac{log_*(\mu_{\eta}k_2)}{\mu_n^2}$.
Recalling that $k+\frac{1}{\mu_{\eta}} \leq \frac{log_*(k_1)}{\mu_{\eta}^2}-\frac{log_*(k_2)}{\mu_{\eta}^2}$,
it's enough to show that $\frac{log_*(\mu_{\eta}k_2)}{\mu_{\eta}^2}-\frac{log_*(k_2)}{\mu_{\eta}^2} \leq \frac{1}{\mu_{\eta}}$.

Case 1: $k_2 \leq \mu_{\eta}$. In this case, it's enough to show
that $log_*(\mu_{\eta}k_2)-log_*(k_2) \leq \mu_{\eta}$, and indeed,
$log_*(\mu_{\eta}k_2)-log_*(k_2) \leq log_*(\mu_{\eta}^2) \leq \mu_{\eta}$.

Case 2: $\mu_{\eta}<k_2$. By the properties of $log_*$, $log_*(k_2) \leq log_*(\mu_{\eta}k_2) \leq log_*(k_2^2) \leq log_*(k_2)+1$,
therefore $\frac{log_*(\mu_{\eta}k_2)}{\mu_{\eta}^2}-\frac{log_*(k_2)}{\mu_{\eta}^2} \leq \frac{1}{\mu_{\eta}}$.

$\square$

\textbf{Definition} \textbf{5}: For $\mathbf n \in \mathbf N$ we define
$\mathbf m=\mathbf{m}_{\mathbf n}^{2}=(\mathbb{Q}_{\mathbf n}^{2},\underset{\sim}{\eta_{\mathbf n}^{2}})$
by:

A) $p\in \mathbb{Q}_{\mathbf n}^{2}$ iff $p$ consists of a trunk $tr(p)\in T_{\mathbf n}$,
a perfect subtree $T_p \subseteq T_{\mathbf n}^{[tr(p)\leq]}$ and a
natural number $n\in [1,lg(tr(p))+1]$ such that $1+ \frac{1}{n} \leq nor_{\eta}(suc_{T_p}(\eta))$
for every $\eta \in T_p^+$.

B) Order: reverse inclusion.

C) $\underset{\sim}{\eta_{\mathbf n}^{2}}=\cup \{tr(p) : p\in \underset{\sim}{G_{\mathbb{Q}_{\mathbf n}^{2}}}\}$.

D) If $p\in \mathbb{Q}_{\mathbf n}^{2}$ we let $nor(p)=min \{n : \eta \in T_p \rightarrow 1+\frac{1}{n}\leq nor_{\eta}(suc_p(\eta))\}$.

\textbf{Claim} \textbf{6:} $\mathbb{Q}_{\mathbf n}^{i} \models ccc$
for $i\in \{0,\frac{1}{2},1,2\}$.

\textbf{Proof}: First we shall prove the claim for $\mathbb{Q}_{\mathbf n}^i$
where $i\in \{0,\frac{1}{2},1\}$. Observe that if $p\in \mathbb{Q}_{\mathbf n}^i$
and $0<k<\omega$, then there is $p\leq q \in \mathbb{Q}_{\mathbf n}^i$
such that $nor_{\eta}(Suc_q(\eta))>k$ for every $\eta \in T_q^+$.
The statement is trivial for $i=1$, so suppose that $i\in \{0,\frac{1}{2}\}$.
In order to prove this fact, let $Y=\{ \eta \in T_p : $for every
$\eta \leq \nu \in T_p$, $nor_{\nu}(Suc_{T_p}(\nu))>k \}$, then
$Y$ is dense in $T_p$ (suppose otherwise, then we can construct
a strictly increasing sequence of memebrs $\eta_i \in T_p$ such that
$nor_{\eta_i}(Suc_{T_p}(\eta_i)) \leq k$, so $\underset{i<\omega}{\cup}\eta_i \in lim(T_p)$
contradicts the definition of $\mathbb{Q}_{\mathbf n}^i$). Now pick
$tr(p) \leq \eta \in Y$, then $q=p^{[\eta \leq]}$ is as required.
It also follows that from this claim that $\mathbb{Q}_{\mathbf n}^{\frac{1}{2}}$
is dense in $\mathbb{Q}_{\mathbf n}^0$.

Now suppose towards contradiction that $\{p_{\alpha} : \alpha<\aleph_1\} \subseteq \mathbb{Q}_{\mathbf n}^i$
is an antichain, for every $\alpha$, there is $p_{\alpha} \leq q_{\alpha}$
such that $nor_{\eta}(Suc_{q_{\alpha}}(\eta))>2$ for every $\eta \in q_{\alpha}$.
For some uncountable $S\subseteq \aleph_1$, $tr(q_{\alpha})=\eta_*$
for every $\alpha \in S$. By the claim below, $q_{\alpha}, q_{\beta}$
are compatible for $\alpha,\beta \in S$, contradicting our assumption.

As for $\mathbb{Q}_{\mathbf n}^2$, given $I=\{p_i : i<\aleph_1\}\subseteq \mathbb{Q}_{\mathbf n}^{2}$
$(\mathbb{Q}_{\mathbf n}^{1})$, the set $\{(tr(p),nor(p)) : p\in I\}$
is countable, hence there is $p_*\in I$ such that for uncountably
many $p_i\in I$ we have $(tr(p_i),nor(p_i))=(tr(p_*),nor(p_*))$.
By the claim below, those $p_i$ are pairwise compatible.

\noindent \begin{flushright}
$\square$
\par\end{flushright}

\textbf{Remark}: The above argument actually gives $\sigma$-linkedness, though this is not used in the paper.

\textbf{Claim} \textbf{7}: 1) $p,q\in \mathbb{Q}_{\mathbf n}^{2}$ are
compatible in $\mathbb{Q}_{\mathbf n}^{2}$ iff $tr(p)\leq tr(q)\in T_p$
or $tr(q)\leq tr(p)\in T_q$.

2) Similarly, $p,q \in \mathbb{Q}_{\mathbf n}^i$ are compatible in
$\mathbb{Q}_{\mathbf n}^i$ for $i\in \{0,\frac{1}{2},1\}$ iff $tr(p) \leq tr(q) \in T_p \vee tr(q) \leq tr(p) \in T_q$.

\textbf{Proof}: In both clauses, the implication $\rightarrow$ is
obvious, we shall prove the other direction.

1) First observe that if $p\in \mathbb{Q}_{\mathbf n}^{2}$ and $\nu \in T_p$,
then $p^{[\nu]}\in \mathbb{Q}_{\mathbf n}^{2}$ and $p\leq p^{[\nu]}$
(where $p^{[\nu]}$ is the set of nodes in $p$ comparable with $\nu$). 

$\square_1$ If $tr(p) \leq tr(q) \in T_p$ then $T_p \cap T_q$ has
arbitrarily long sequences.

Proof: Let $\eta=tr(q)$, then by the definition of the norm and $\mathbb{Q}_{\mathbf n}^{2}$,
$\frac{1}{2}<\frac{|suc_{T_p}(\eta)|}{|suc_{T_{\mathbf n}}(\eta)|},\frac{|suc_{T_q}{\eta}|}{|suc_{T_{\mathbf n}}(\eta)|}$.
Hence there is $\nu \in suc_{T_p}(\eta)\cap suc_{T_q}(\eta)$. Repeating
the same argument, we get sequences in $T_p \cap T_q$ of length $n$
for every $n$ large enough.

$\square_2$ Claim: If $tr(p_1)=tr(p_2)=\eta$, $p_1,p_2 \in \mathbb{Q}_{\mathbf n}^{2}$,
$min \{ nor(p_1), nor(p_2) \} \leq h$ and $h<lg(\eta)$, then $p_1$
and $p_2$ are compatible.

Proof: For every $\nu \in T_{p_1}\cap T_{p_2}$, by the co-bigness
property, $min\{nor_{\nu}(suc_{p_1}(\nu)), suc_{p_2}(\nu)\}-\frac{1}{\mu_{\nu}}\leq nor(suc_{p_1}(\nu) \cap suc_{p_2}(\nu))$.
By the definition of $nor(p_i)$ (recalling that $lg(\eta)^2 \leq \mu_{\eta}$),
$1+\frac{1}{h+1} \leq (1+\frac{1}{h+1})+(\frac{1}{(h+1)^2}-\frac{1}{\mu_{\eta}}) \leq (1+\frac{1}{h+1})+(\frac{1}{h}-\frac{1}{h+1}-\frac{1}{\mu_{\nu}})=1+\frac{1}{h}-\frac{1}{\mu_{\nu}} \leq min\{nor(p_1),nor(p_2)\}-\frac{1}{\mu_{\nu}} \leq min\{nor_{\nu}(suc_{p_1}(\nu)), suc_{p_2}(\nu)\}-\frac{1}{\mu_{\nu}}$.
Therefore $1+\frac{1}{h+1}\leq nor(suc_{p_1}(\nu) \cap suc_{p_2}(\nu))$,
so $p_1 \cap p_2$ is as required. Hence:

$\square_3$ $p$ and $q$ are compatible.

Proof: Suppose WLOG that $tr(p)\leq tr(q) \in T_p$ and pick $h$
such that $1+\frac{1}{h} \leq nor(p),nor(q)$. By $\square_1$, there
is $\eta \in T_p \cap T_q$ such that $h<\lg(\eta)$. Now $p\leq p^{[\eta]}, q\leq q^{[\eta]}$
and $(p^{[\eta]},q^{[\eta]})$ satisfy the assumptions of $\square_2$,
therefore they're compatible and so are $p$ and $q$.

The proof is similar if $tr(q)\leq tr(p)\in T_q$. The implication
in the other direction is easy.

2) The proof is similar. First observe that if $\eta \in lim(T_p) \cap lim(T_q)$,
then $lim(nor_{\eta \restriction l}(suc_{T_p}(\eta \restriction l)) : l<\omega)=\infty=lim(nor_{\eta \restriction l}(suc_{T_q}(\eta \restriction l)) : l<\omega)$,
so by the co-bigness property (definition 2(f)), $lim(nor_{\eta \restriction l}(suc_{T_{p\cap q}}(\eta \restriction l)) : l<\omega)=\infty$.
Now let $\nu=tr(q) \in T_p \cap T_q$, as $2-\frac{1}{\mu_{tr(p)}} \leq nor(p),nor(q)$,
it follows from the co-bigness property and definition 2(g) that $\nu \leq \eta \in T_p \cap T_q \rightarrow 2<|Suc_{p\cap q}(\eta)|$,
so $p\cap q$ is a perfect tree. It's easy to see that there exists
$\eta \in p\cap q$ such that $nor_{\nu}(Suc_{p\cap q}(\nu))>2$ for
every $\eta \leq \nu \in p\cap q$ (otherwise, we can repeart the
argument in the proof of claim 6, and get a branch through $p\cap q$
along which the norm doesn't tend to infinity). Therefore, $p^{[\leq \eta]} \cap q^{[\leq \eta]} \in \mathbb{Q}_{\mathbf n}^i$
$(i\in \{0,\frac{1}{2}\})$ is a common upper bound. Finally, note
that if $i=1$, then for every $n<\omega$ there exist $k^p(n+1),k^q(n+1)$
as in definition 3.3. By the co-bigness property, for every $\eta \in T_p \cap T_q$
of length $>max\{k^p(n+1),k^q(n+1)\}$, $n\leq nor_{\eta}(Suc_{p\cap q}(\eta))$.
Therefore, the common upper bound is in $\mathbb{Q}_{\mathbf n}^1$
as well. 

\noindent \begin{flushright}
$\square$
\par\end{flushright}

\textbf{Claim} \textbf{8}: Let $I\subseteq \mathbb{Q}_{\mathbf n}^{2}$
be an antichain and $A=\cup \{T_q^+ : q\in I\}\subseteq T_{\mathbf n}$.
The following conditions are equivalent:

(a) $I$ is a maximal antichain.

(b) If $\eta \in T_{\mathbf n}$ and $0<n<\omega$ then there is no
$p\in \mathbb{Q}_{\mathbf n}^{2}$ such that:

$(\alpha)$ $tr(p)=\eta$. 

$(\beta)$ $nor(p)=n$.

$(\gamma)$ $p$ is incompatible with every $q\in I$. 

(c) Like (b), but replcaing $(\gamma)$ by

$(\gamma)'$ $T_p^+\cap A=\emptyset$.

(d) Like (b), but replcaing $(\gamma)$ by

$(\gamma)''$ For every $m>n$ $T_p^+\cap A$ is disjoint to $\{\nu \in T_{\mathbf n} : lg(\nu)\leq m\}$. 

(e) If $\eta \in T_{\mathbf n}$ and $n<\omega$ then for some $m>n$
there is no set $T$ such that:

$(\alpha)$ $T\subseteq T_{\mathbf n}$.

$(\beta)$ $\eta \in T$.

$(\gamma)$ If $\nu \in T^+$ then $\eta \leq \nu$ and $lg(\nu)\leq m$.

$(\delta)$ If $\eta \leq \nu_1 \leq \nu_2$ and $\nu_2\in T$ then
$\nu_1 \in T$.

$(\epsilon)$ $T\cap A=\emptyset$.

$(\zeta)$ If $\nu \in T$ and $lg(\nu)<m$ then $1+\frac{1}{n} \leq nor_{\nu}(suc_T(\nu))$.

\textbf{Proof}: $\neg (a) \rightarrow \neg (b):$ If $p$ is incompatible
with every $q\in I$ then $(p,tr(p),nor(p))$ is a counterexample
to (b).

$\neg (b)\rightarrow \neg (c):$ If $(p,tr(p),nor(p))$ is a counterexample
to (b), then it is a counterexample to (c) by the characterisation
of compatibility in $\mathbb{Q}_{\mathbf n}^{2}$ in claim 7.

$\neg (c) \rightarrow \neg (d):$ Obvious.

$\neg (d) \rightarrow \neg (e)$: Let $T=T_p$ with $p$ being a counter
example to $(d)$ and let $\eta=tr(p), n$ witness $\neg(d)$. We
shall check that for every $m>n$, $\{ \nu : tr(p) \leq \nu \in T \wedge lg(\nu) \leq m\}$
satisfies $(\alpha)-(\zeta)$ if $(e)$.

$\neg (e) \rightarrow \neg (a):$ If $(\eta, n)$ is a counterexample,
then for every $m$ there is $T_m$ satisfying $(\alpha)-(\zeta)$
of clause $(e)$. Let $D$ be a non-principal ultrafilter on $\omega$
and define $T:=\{\nu \in T_{\mathbf n} : \nu \leq \eta$ or $\{m : m>n, \nu \in T_m\}\in D\}$.
It remains to show that $T\in \mathbb{Q}_{\mathbf n}^{2}$ (as $T^+$
is disjoint to $A$, it follows that $I$ is not a maximal antichain).
The proof is similar to claim 12. 

\noindent \begin{flushright}
$\square$
\par\end{flushright}

\textbf{Claim} \textbf{9}: Let $\mathbf n \in \mathbf N$.

A) The sets $\mathbb{Q}_{\mathbf n}^1$ and $\mathbb{Q}_{\mathbf n}^2$
are Borel, the sets $\mathbb{Q}_{\mathbf n}^0$ and $\mathbb{Q}_{\mathbf n}^{\frac{1}{2}}$
are $\Pi_1^1$.

B) The relation $\leq_{\mathbb{Q}_{\mathbf n}^i}$ is Borel for $i\in \{0,\frac{1}{2},1,2\}$.

C) The incompatibility relation in $\mathbb{Q}_{\mathbf n}^i$ is Borel
for $i\in \{0,\frac{1}{2},1,2\}$.

\textbf{Proof:}

\textbf{A. The sets $\mathbb{Q}_{\mathbf n}^1$ and $\mathbb{Q}_{\mathbf n}^2$
are Borel:} We shall first prove the claim for $\mathbb{Q}_{\mathbf n}^1$.
Consider $T_{\mathbf n}$ as a subset of $H(\aleph_0)$. By definition,
if $p\in \mathbb{Q}_{\mathbf n}^{1}$ then $T_p \subseteq T_{\mathbf n}\subseteq H(\aleph_0)$.
Hence $S:=\{p\subseteq H(\aleph_0) : p$ is a perfect subtree of $T_{\mathbf n} \}\subseteq P(H(\aleph_0))$
is a Borel subset of $P(H(\aleph_0))$. For every $n,k<\omega$ define
$S_{n,k}^{1}=\{p\in S : lg(tr(p))<k$ and if $\rho \in T_p$ and $k\leq lg(\rho)$
then $n\leq nor_{\rho}(suc_p(\rho))\}$. Each $S_{n,k}^{1}$ is closed,
hence $S \cap (\underset{n}{\cap}\underset{k}{\cup}S_{n,k}^{1})$
is Borel, so it's enough to show that $p\in \mathbb{Q}_{\mathbf n}^{1}$
iff $p\in S \cap (\underset{n}{\cap}\underset{k}{\cup}S_{n,k}^{1})$
and $2-\frac{1}{\mu_{tr(p)}} \leq nor(p)$, which follows directly
from the definition of $\mathbb{Q}_{\mathbf n}^1$. 

In the case of $\mathbb{Q}_{\mathbf n}^{2}$, we replace $\underset{n}{\cap}\underset{k}{\cup}S_{n,k}^{1}$
with $\underset{n,k}{\cup}S_{n,k}^{2}$ where $S_{n,k}^{2}=\{p\in S : lg(tr(p))=n \wedge nor(p)=k\}$.
Each $S_{n,k}^{2}$ is Borel and since {}``being a perfect subtree''
is Borel, $\mathbb{Q}_{\mathbf n}^{2}$ is Borel. 

\textbf{The sets $\mathbb{Q}_{\mathbf n}^0$ and $\mathbb{Q}_{\mathbf n}^{\frac{1}{2}}$
are $\Pi_1^1$: }The demand {}``$\underset{n<\omega}{lim}(nor_{\eta \restriction n}(Suc_p(\eta \restriction n)))=\infty$
for every $\eta \in lim(T_p)$'' is $\Pi_1^1$,and it's easy to see
that $\{p\in S : tr(p) \leq \eta \in T_p \rightarrow nor_{\eta}(Suc_{T_p}(\eta))>2\}$
is Borel.

\textbf{B. The relation $\leq_{\mathbb{Q}_{\mathbf n}^i}$ is Borel
for $i\in \{0,\frac{1}{2},1,2\}$: }For $i\in \{0,\frac{1}{2},1,2\}$,
the relation $\leq_{\mathbb{Q}_{\mathbf n}^{i}}$ is simply the reverse
inclusion relation restricted to $\mathbb{Q}_{\mathbf n}^{i}$, hence
it is Borel.

\textbf{C. The incompatibility relation in $\mathbb{Q}_{\mathbf n}^i$
is Borel for $\{0,\frac{1}{2},1,2\}$: }The incompatibility relation
is Borel by claim 7.

\noindent \begin{flushright}
$\square$
\par\end{flushright}

\textbf{Claim} \textbf{10}: A) Assume that $p_l \in \mathbb{Q}_{\mathbf n}^{i}$
$(l<n)$ where $i\in \{0,1\}$, $\underset{l<n}{\wedge}tr(p_l)=\rho$,
$n\leq lg(\rho)$ and for every $\eta \in p_l^+$ we have $2\leq k+1\leq nor_{\eta}(suc_{p_l}(\eta))$,
then $\{p_l : l<n\}$ have a common upper bound $p$ such that $tr(p)=\rho$
and $k\leq nor_{\eta}(suc_p(\eta))$ for every $\eta \in T_p^+$.

B) Assume that $p_l \in \mathbb{Q}_{\mathbf n}^{2}$ $(l<n)$, $\underset{l<n}{\wedge}tr(p_l)=\rho$,
$n\leq lg(\rho)$ and for every $\eta \in p_l^+$ $(l<n)$ we have
$1+\frac{1}{k}\leq nor_{\eta}(suc_{p_l}(\eta))$. In addition, assume
that $k\leq lg(\rho)$ and $k(k+1) \leq \mu_{\eta}$ for every $\eta \in p_l^+$ $(l<n)$,
then $\{p_l : l<n\}$ have a common upper bound $p$ such that $tr(p)=\rho$
and $1+\frac{1}{k+1}\leq nor_{\eta}(suc_p(\eta))$.

\textbf{Proof}: A) Suppose first that $i=0$. Let $p=\underset{l<n}{\cap}p_l$,
then $p\subseteq T_{\mathbf n}^{[\rho \leq]}$ is a subtree conatining
$\rho$. If $\nu \in p$ then $\nu \in p_l$ for every $l<n$, hence
$Suc_p(\nu)=\underset{l<n}{\cap}Suc_{p_l}(\nu)$. As $n\leq lg(\rho) \leq \mu_{\eta}$
for every $\rho \leq \eta \in p$, it follows from the properties
of the norm in the definition of $\mathbf n \in \mathbf N$ that $k\leq nor_{\eta}(Suc_p(\eta))$.
Therefore, $T_p$ is a prefect tree, and similarly to the proof of
claim 7, it follows that the norm along infinite branches tends to
infinity, hence $p\in \mathbb{Q}_{\mathbf n}^0$. Suppose now that $i=1$.
The above arguments are still valid, and in addition, similarly to
the argument on $\mathbb{Q}_{\mathbf n}^1$ in th proof of claim 7(2),
it's easy to see that by the co-bigness property, $p\in \mathbb{Q}_{\mathbf n}^1$.

Remark: Note that as $2\leq k+1$, it follows from the above arguments
that $2-\frac{1}{\mu_{tr(p)}} \leq nor_{\eta}(Suc_{T_p}(\eta))$ for
every $tr(p)\leq \eta \in  T_p$. In fact, $k+1-\frac{1}{\mu_{\rho}} \leq nor_{\eta}(Suc_{T_p}(\eta))$,
therefore, if $2<k+1-\frac{1}{\mu_{\rho}}$ then we also get the claim
for $i=\frac{1}{2}$.

B) The proof is similar, the only difference is that now we have to
prove the following assertion:

$(*)$ If $b_l\subseteq suc_{T_{\mathbf n}}(\eta)$ for $l<n\leq \mu_{\eta}$,
$\underset{l<n}{\wedge} 1+\frac{1}{k} \leq nor_{\eta}(b_l)$ and $b=\underset{l<n}{\cap}b_l$
then $1+\frac{1}{k+1}\leq nor_{\eta}(b)$.

The assertion follows from the co-bigness property (definition $2(f)$,
with $b_i$ and $1+\frac{1}{k}-\frac{1}{\mu_{\eta}}$ here standing
for $a_i$ and $k$ there). 

\noindent \begin{flushright}
$\square$
\par\end{flushright}

\textbf{Claim} \textbf{11}: Let $\mathbf n \in \mathbf N$. $"\{p_n : n<\omega\}$
is a maximal antichain'' is Borel for $\{p_n : n<\omega\}\subseteq \mathbb{Q}_{\mathbf n}^{2}$.

\textbf{Proof}: By claim 8.

\noindent \begin{flushright}
$\square$
\par\end{flushright}

\textbf{Claim 12}: Assume $\{p_n : n<\omega\} \subseteq \mathbb{Q}_{\mathbf n}^{2}$,
$\underset{n}{\wedge}tr(p_n)=\eta$ and $\underset{n}{\wedge}nor(p_n)=k$.
Then there is $p_* \in \mathbb{Q}_{\mathbf n}^{2}$ such that: 

(a) $tr(p_*)=\eta$, $nor(p_*)=k$.

(b) $p_* \Vdash_{\mathbb{Q}_{\mathbf n}^{2}} "(\exists^{\infty}n)(p_n \in G_{\mathbb{Q}_{\mathbf n}^{2}})"$.

\textbf{Proof}: Let $D$ be a uniform ultrafilter on $\omega$ and
define $T_{p_*}:=\{\nu \in T_{\mathbf n} : \{n : \nu \in p_n\}\in D\}$.
If $\nu \in T_{p_*}$, then for some $n$, $\nu \in T_{p_n}\subseteq T_{\mathbf n}^{[\eta \leq]}$
(recalling that $tr(p_n)=\eta$), hence $T_{p_*}\subseteq T_{\mathbf n}^{[\eta \leq]}$.
Obviously, $l\leq lg(\eta) \rightarrow \eta \restriction l \in T_{p_*}$
as $\eta=tr(p_n)\in p_n$ for every $n$. 

$(*)_1$ If $\eta \triangleleft \nu \triangleleft  \rho$ and $\rho \in T_{p_*}$,
then $\nu  \in T_{p_*}$. 

Why? Define $A_{\rho}=\{n : \rho \in p_n\}$ and define $A_{\nu}$
similarly. $A_{\rho} \in D$ by the definition of $T_{p_*}$. Obviously
$A_{\rho}\subseteq A_{\nu}$, hence $A_{\nu} \in D$ and $\nu \in T_{p_*}$.

$(*)_2$ If $\nu \in T_{p_*}$ then $1+\frac{1}{k}\leq nor_{\nu}(suc_{p_*}(\nu))$.

Why? Define $A_{\nu}$ as above, so $A_{\nu}\in D$. Let $(b_l : l<l(*))$
list $\{suc_{p_n}(\nu) : n<\omega\}$. As $\{suc_{p_n}(\nu) : n<\omega\} \subseteq P(suc_{T_{\mathbf n}}(\nu))$,
we have $l(*)\leq 2^{|suc_{T_{\mathbf n}}(\nu)|}=2^{\lambda_{\nu}}<\aleph_0$.
For $l<l(*)$ let $A_{\nu,l}:=\{n \in A_{\nu} : suc_{p_n}(\nu)=b_l\}$.
Obviously this is a finite partition of $A_{\nu}$, hence there is
exactly one $m<l(*)$ such that $A_{\nu,m}\in D$ and therefore $b_m \subseteq suc_{p_*}(\nu)$
and actually $b_m=suc_{p_*}(\nu)$ (if $\eta \in suc_{p_*}(\nu)$
is witnessed by $X\in D$, then $X\cap A_{\nu,m}$ is a witness for
$\eta \in b_m$). Therefore $nor_{\nu}(b_m)= nor_{\nu}(suc_{p_*}(\nu))$
and for some $n$ we have $1+\frac{1}{k}=1+\frac{1}{nor(p_n)}\leq nor_{\nu}(suc_{p_n}(\nu))= nor_{\nu}(suc_{p_*}(\nu))$. 

It follows from the above arguments that $p_* \in \mathbb{Q}_{\mathbf n}^2$.

We shall now prove that

$(*)_3$$p_* \Vdash_{\mathbb{Q}_{\mathbf n}^{2}} "(\exists^{\infty}n)(p_n \in G_{\mathbb{Q}_{\mathbf n}^{2}})"$. 

Why? Suppose that $p_*\leq q$, then $tr(q)\in T_{p_*}$. By the definition
of $p_*$, $\{n : tr(q)\in p_n\}\in D$. For every such $p_n$, $\eta=tr(p_n) \leq tr(q) \in T_{p_n}$,
so $p_n$ is compatible with $q$ and hence with $p_*$.

\noindent \begin{flushright}
$\square$
\par\end{flushright}

\textbf{Claim 12': }For $\iota \in \{0,\frac{1}{2},1,2\}$, $\underset{\sim}{\eta_{\mathbf n}^{\iota}}$
is a generic for $\mathbb{Q}_{\mathbf n}^{\iota}$, i.e. $\Vdash_{\mathbb{Q}_{\mathbf n}^{\iota}} "V[\underset{\sim}{G_{\mathbb{Q}_{\mathbf n}^{\iota}}}]=V[\underset{\sim}{\eta_{\mathbf n}^{\iota}}]"$. 

\textbf{Proof: }Easy.

\noindent \begin{flushright}
$\square$
\par\end{flushright}

\textbf{\large 3. The iteration}{\large \par}

In this section we shall describe our iteration. Although our definition
will be general and will follow the technique of iteration along templates
as described in {[}Sh700{]}, we will eventually use a simple private
case of the general construction (see also [Br] and [Me2]). In our case, we'll have a non-wellfounded
linear order $L$, and the forcing will be the union of finite-length
iterations along subsets of $L$. Dealing with FS-iterations of Suslin
forcing will guarantee that the union is well-behaved.

\textbf{Iteration parameters}

The purpose of Definitions 12 and 13 is to show how our construction fits as a special case in the broader context of the second author's general method of iterations along templates. However, we shall only use the private case of Definition 13(A), and so a reader who only wants to understand the main results in this paper may focus on Definition 13(A).

\textbf{Definition} \textbf{12}: Let $\mathbf Q$ be the class of $\mathbf q$
(iteration parameters) consisting of:

a. A partial order $L_{\mathbf q}=L[\mathbf q]$.

b. $\bar{u}_0=(u_t^0 : t\in L_{\mathbf q})$ such that $u_t^0 \subseteq L_{<t}$
for each $t\in L_{\mathbf q}$ (and $u_t^0$ is well-ordered by (d)).
In the main case $|u_t^0|\leq \aleph_0$ (in our application, $u_t^0$
is actually empty).

c. $\mathbf{I}=(\mathbf{I}_t : t\in L_{\mathbf q})$ such that each $\mathbf{I}_t$
is an ideal on $L_{<t}$ and $u_t^0\in \mathbf{I}_t$. In the main case
here, $\mathbf{I}_t=\{u\subseteq L_{<t} : u $ is finite$\}$.

d. $\mathbf L$ is a directed family of well-founded subsets of $L_{\mathbf q}$
closed under initial segments such that $\underset{L\in \mathbf{L}}{\cup}L=L_{\mathbf q}$
and $t\in L \rightarrow u_t^0 \subseteq L$ (for $L\in \mathbf{L}$).

e. $(\mathbf{m}_t : t\in L_{\mathbf q})$ is a sequence such that each
$\mathbf{m}_t$ is a definition of a Suslin ccc forcing notion $\mathbb{Q}_{\mathbf{m}_t}^i$
with a generic $\underset{\sim}{\eta_{\mathbf{m}_t}}$ (depending on
a formula using $\mathbf{B}_t(...,\eta_s,...)_{s\in u_t^0}$, see f+g
and definition 13).

f. Actually, $\mathbf{m}_t=\mathbf{m}_{t,\underset{\sim}{\nu}_t}$ where
$\underset{\sim}{\nu_t}=\mathbf{B}_t(\bar{\eta} \restriction u_t^0)$
is a name of a real and $\mathbf{B}_t$ is a Borel function (see definition
13(E) below), i.e. $\mathbf{m}_t$ is computed from the parameter $\underset{\sim}{\nu_t} \in \omega^{\omega}$.

g. For every $t\in L_{\mathbf q}$, $\mathbf{B}_t : \underset{i\in u_t^0}{\Pi}\omega^{\omega} \rightarrow \omega^{\omega}$
is an absolute Borel function.

h. For a linear order $L$, let $L^+:=L \cup \{\infty\}$ which is
obtained by adding an element above all elements of $L$.

\textbf{The iteration}

\textbf{Definition} \textbf{and claim 13}: For $i\in \{1,2\}$, $\mathbf{q} \in \mathbf Q$
and $L\in \mathbf L$ we shall define the FS iteration $\bar{\mathbb{Q}}_L=(\mathbb{P}_t^L, \underset{\sim}{\mathbb{Q}_t^L} : t\in L^+)$
with limit $\mathbb{P}_L$ and the $\mathbb{P}_t^L=\mathbb{P}_{L,<t}$-names
$\underset{\sim}{\eta_t}, \underset{\sim}{\nu_t}$ by induction on
$dp(L)$ (where $dp(L)$ is the depth of $L$, recalling that $L$
is well-founded) such that:

A. a) $\mathbb{P}_L$ is a forcing notion.

b) $\underset{\sim}{\eta_t}$ is a $\mathbb{P}_L$ name when $u_t^0 \cup \{t\}\subseteq L\in \mathbf L$
(so we use a maximal antichain from $\mathbb{P}_L$, moreover, from
$\mathbb{P}_{L_1}$ for every $L_1 \in \mathbf L$ which is $\subseteq L$).

c) $\underset{\sim}{\nu_t}$ is a $\mathbb{P}_L$ name when $u_t^0  \subseteq L\in \mathbf{L}$.

d) If $L_1,L_2 \in \mathbf{L}$ are linearly ordered, $L_1 \subseteq L_2$
and each $\mathbf{I}_t$ has the form $\{ L\subseteq L_{<t} : L$ is
well-ordered$\}$, then $\mathbb{P}_{L_1}\lessdot \mathbb{P}_{L_2}$.

B. $p\in \mathbb{P}_t^L$ iff

a. $Dom(p)\subseteq L_{<t}$ is finite.

b. If $s\in Dom(p)$ then for some $u\in \mathbf{I}_s \cap \mathcal{P}(L_{<s})$
and a Borel function $\mathbf B$, $p(s)=\mathbf{B}(...,\underset{\sim}{\eta_r},...)_{r\in u}$
and $\Vdash_{\mathbb{P}_s^L} "p(s)\in \mathbb{Q}_{\mathbf{m}_s}^{i}"$.

C. $\underset{\sim}{\mathbb{Q}_{t}^{L}}$ is the $\mathbb{P}_{t}^{L}$-name
of $\mathbb{Q}_{\mathbf{m}_t}^{i}$ using the parameter $\underset{\sim}{\nu_t}$.

D. $\bar{\eta}=(\underset{\sim}{\eta_t} : t\in L_{\mathbf q})$. Each
$\underset{\sim}{\eta_t}$ is defined as the generic of $\mathbb{Q}_t^L$
(by a maximal antichain of $\mathbb{P}_{L}$ whenever $L\in \mathbf{L}$
and $u_t^0 \subseteq L \subseteq L_{<t}$), meaning: $t\in L \in \mathbf L \rightarrow \Vdash "\underset{\sim}{\eta_t}$
is a generic for $\underset{\sim}{\mathbb{Q}_t}"$ defined as usual.

E. $\bar{\nu}=(\underset{\sim}{\nu_t} : t\in L_{\mathbf q})$ such that
for each $t\in L_{\mathbf q}$, $\mathbf{B}_t$ is a Borel function and
$\underset{\sim}{\nu_t}=\mathbf{B}_t(\bar{\eta}\restriction u_t^0)$.

F. The order on $\mathbb{P}_L$ is defined naturally.

\textbf{Proof: }Should be clear.

\noindent \begin{flushright}
$\square$
\par\end{flushright}

\textbf{13(A) A special case of the general construction}

Of special interest here is the case where $\mathbf q \in \mathbf Q$
satisfies:

a. $L_{\mathbf q}$ is a dense linear order, $\mathbf{I}_t=[L_{<t}]^{<\aleph_0}$
for each $t\in L_{\mathbf q}$ and $\mathbf{L}=[L_{\mathbf q}]^{<\aleph_0}$.

b. $\mathbf{m}_t$ is a definition of $\mathbb{Q}_{\mathbf{n}_t}^{i}$
where $i\in \{1,2\}$ (hence a Suslin c.c.c. forcing), not using a
name of the form $\underset{\sim}{\nu_t}$.

c. $\mathbf{m}_t \in V$ and $u_t^0=\emptyset$ for every $t\in L_{\mathbf q}$.

\textbf{13(B) }We shall denote the collection of $\mathbf{q} \in \mathbf{Q}$
as above by $\mathbf{Q}_{sp}$.

\textbf{13(C) Hypothesis: }From now on we assume that $\mathbf q \in \mathbf Q$
satisfies the requirements of $13(A)$.

\textbf{Definition/Observation} \textbf{14}: Let $\mathbf q \in \mathbf Q$.

1. $\{\mathbb{P}_{J} : J\subseteq L_{\mathbf q}$ is finite$\}$ is
a $\lessdot$-directed set of forcing notions.

2. For $J\subseteq L_{\mathbf q}$, let $\mathbb{P}_J=\cup \{ \mathbb{P}_{J'} : J' \subseteq J$
is finite$\}$ and $\mathbb{P}_{\mathbf q}=\mathbb{P}_{L_{\mathbf q}}$.

\textbf{Proof}: (1) follows by {[}JuSh292{]}.

\noindent \begin{flushright}
$\square$
\par\end{flushright}

\textbf{Claim} \textbf{15}: 1) For every $J_1\subseteq J_2\subseteq L_{\mathbf q}$,
$\mathbb{P}_{J_1}\lessdot \mathbb{P}_{J_2}$.

2) If $J\subseteq L_{\mathbf q}$ then $\mathbb{P}_{J}=\mathbb{P}_{\mathbf q,J}=\cup \{\mathbb{P}_I : I\subseteq J$
is finite$\} \lessdot \mathbb{P}_{\mathbf q}$.

\textbf{Proof: }1) Case 1: $|J_2|<\aleph_0$. Easy by {[}JuSh292{]}.

Case 2: $J_2$ is inifinite. Let $q\in \mathbb{P}_{J_2}$, then for
some finite $J_2^* \subseteq J_2$, $q\in \mathbb{P}_{J_2^*}$. Let
$J_1^*=J_1 \cap J_2^*$. As $\mathbb{P}_{J_1^*} \lessdot \mathbb{P}_{J_2^*}$
by observation 14(1), there is $p\in \mathbb{P}_{J_1^*}$ such that
$p\leq p' \in \mathbb{P}_{J_1^*} \rightarrow p'$ and $q$ are compatible.
It suffices to prove that if $J_1' \subseteq J_1$ is finite and $J_1^* \subseteq J_1'$,
then $p\leq p' \in \mathbb{P}_{J_1'} \rightarrow p'$ and $q$ are
compatible in $\mathbb{P}_{J_2^* \cup J_1'}$ (as if $p\leq p' \in \mathbb{P}_{J_1}$,
then $p' \in \mathbb{P}_{J_1'}$ where $J_1'=J_1^* \cup Dom(p')$).
We prove this by induction on $sup\{|L_{<t} \cap J_1^*| : t\in J_1' \setminus J_1^*\}$
as in {[}JuSh292{]}.

2) By (1).

\textbf{Observation} \textbf{16}: Suppose that $\mathbf q \in \mathbf Q$,
$J\in \mathbf{L}$ is finite and $p_1,p_2 \in \mathbb{P}_{J}$. If $tr(p_1(t))=tr(p_2(t))$
for every $t\in Dom(p_1)\cap Dom(p_2)$, then $p_1$ and $p_2$ are
compatible.

\textbf{Proof}: By induction on $|J|$. The induction step is a corollary
of the compatibility condition for $\mathbb{Q}_{\mathbf n}^{2}$ (see
claim 7).

\noindent \begin{flushright}
$\square$
\par\end{flushright}

\textbf{Claim 17: }For $\mathbf q \in \mathbf Q$, $\mathbb{P}_{\mathbf q} \models ccc$.

\textbf{Proof: }Suppose that $\{p_{\alpha} : \alpha<\aleph_1\} \subseteq \mathbb{P}_{\mathbf q}$.
For each $\alpha<\aleph_1$ there is a finite $J_{\alpha} \subseteq L_{\mathbf q}$
such that $p_{\alpha} \in \mathbb{P}_{J_{\alpha}}$. Hence there is
$n_*\in \mathbb{N}$ such that $|\{p_{\alpha} : |J_{\alpha}|=n_*\}|=\aleph_1$.
For each $\alpha$ denote $J_{\alpha}=\{t_{\alpha,0}<...<t_{\alpha,n_{\alpha-1}}\}$,
by cardinallity arguments i.e. the $\Delta$-system lemma, WLOG there
is $u\subseteq n_*$ such that $t_{\alpha,l}=t_l$ for every $\alpha<\aleph_1$
and $(t_{\alpha,l} : l\in n_* \setminus u, \alpha<\aleph_1)$ is without
repetitions. As every condition $p_{\alpha} \in \mathbb{P}_{J_{\alpha}}$
belongs to an iteration along $J_{\alpha}$ in the usual sense, there
is $p_{\alpha} \leq p_{\alpha}' \in \mathbb{P}_{J_{\alpha}}$ such
that $tr(p_{\alpha}'(t))$ is an object for every $t\in J_{\alpha}$
(so $J_{\alpha}=Dom(p_{\alpha}')$). Given $l\in u$ there are countably
many possible values for $tr(p_{\alpha}(t_l))$, hence there is a
set $I=\{p_{{\alpha}_i} : i<i(*)\}\subseteq \{p_{\alpha} : \alpha<\aleph_1\}$
of cardinality $\aleph_1$ such that $tr(p_{{\alpha}_i}(t_l))$ is
constant for all $i<i(*)$. If $i<j<i(*)$, then $J_{i,j}:=J_{\alpha_i} \cup J_{\alpha_j} \subseteq L_{\mathbf q}$
is finite, $p_{\alpha_i} \in  \mathbb{P}_{J_{\alpha_i}} \lessdot \mathbb{P}_{J_{i,j}}$
and $p_{\alpha_j} \in  \mathbb{P}_{J_{\alpha_j}} \lessdot \mathbb{P}_{J_{i,j}}$,
so $p_{\alpha_i}$ and $p_{\alpha_j}$ are compatible in $\mathbb{P}_{J_{i,j}}$
(hence in $\mathbb{P}_{\mathbf q}$) by observation 16.

\noindent \begin{flushright}
$\square$
\par\end{flushright}

\textbf{\large 4. The ideals derived from a forcing notion $\mathbb Q$}{\large \par}

We shall now define the ideals derived from a Suslin forcing notion
$\mathbb{Q}$ and a name $\underset{\sim}{\eta}$ of a real.

\textbf{Definition} \textbf{18}: 1. Let $\mathbb Q$ be a forcing
notion such that each $p\in \mathbb Q$ is a perfect subtree of $\omega^{<\omega}$,
$p\leq_{\mathbb Q} q$ iff $q\subseteq p$ and the generic real is
given by the union of trunks of conditions that belong to the generic
set, that is $\underset{\sim}{\eta}=\underset{p\in \underset{\sim}{G}}{\cup} tr(p)$
and $\Vdash_{\mathbb Q} "\underset{\sim}{\eta} \in \omega^{\omega}"$.

Let $\aleph_0 \leq \kappa$, the ideal $I_{\mathbb Q,\kappa}^0$ will
be defined as the closure under unions of size $\leq \kappa$ of sets
of the form $\{X\subseteq \omega^{\omega} : (\forall p\in \mathbb Q)(\exists q\geq p)(lim(q) \cap X=\emptyset)\}$.\footnote{The above definition has the following variant in the literature, which will not be used in this paper: Let $\mathbf m=(\mathbb Q,\underset{\sim}{\eta})$ where $\underset{\sim}{\eta}$
is a $\mathbb Q$-name of a real, the ideal $I_{\mathbf m,\kappa}^1$
for $\aleph_0 \leq \kappa$ will be defined as follows:
\\
$A\in I_{\mathbf m,\kappa}^1$ iff there exists $X\subseteq \kappa$
such that $A\cap \{ \underset{\sim}{\eta}[G] : G\subseteq \mathbb{Q}^{L[X]}$
is generic over $L[X]\}=\emptyset$.}

2. For $\mathbb{Q}$ as above, we let $I_{\mathbb Q, <\aleph_0}$ be the set $\{X\subseteq \omega^{\omega} : (\forall p\in \mathbb Q)(\exists q\geq p)(lim(q) \cap X=\emptyset)\}$.

3. For $(\mathbb Q, \underset{\sim}{\eta})$ and $\kappa$ as in (1), we shall denote $I_{\mathbb Q,\kappa}^0$
by $I_{\mathbb Q,\kappa}$.

4. Let $I$ be an ideal on the reals, a set of reals $X$ is called
$I$-measurable if there exists a Borel set $B$ such that $X\Delta B \in I$.

5. A set of reals $X$ will be called $(\mathbb Q,\kappa)$-measurable
if it is $I_{\mathbb Q,\kappa}$-measurable.

6. Given a model $V$ of $ZF$, we say that $(\mathbb Q,\kappa)$-measurability
holds in $V$ if every set of reals in $V$ is $(\mathbb Q,\kappa)$-measurable
and $I_{\mathbb Q,\kappa}$ is a non-trivial ideal.

Remark: In {[}F1424{]} we shall further investigate the above ideals.

\textbf{\large 5. Cohen reals}{\large \par}

An important feature of $\mathbb{Q}_{\mathbf n}^{\iota}$ is the fact
that it adds a Cohen real. This fact will be later used to show that
$\mathbb{Q}_{\mathbf n}^{\iota}$ can turn the ground model reals into
a null set with respect to the relevant ideal.

\textbf{Claim} \textbf{19}: Forcing with $\mathbb{Q}_{\mathbf n}^{\iota}$
$(i\in \{0,\frac{1}{2},1,2\})$ adds a Cohen real.

\textbf{Proof}: For every $\eta \in T_{\mathbf n}$ let $g_{\eta}: suc_{T_{\mathbf n}}(\eta)\rightarrow \{0,1\}$
be a function such that $|g_{\eta}^{-1}\{l\}|>\frac{\lambda_{\eta}}{2}-1$
$(l=0,1)$ (recall that $\lambda_{\eta}=|suc_{T_{\mathbf n}}(\eta)|$).
Define a $\mathbb{Q}_{\mathbf n}^{\iota}$-name $\underset{\sim}{\nu}$
by $\underset{\sim}{\nu}(n)=g_{\underset{\sim}{\eta_{\mathbf n}^{\iota}} \restriction n}(\underset{\sim}{\eta_{\mathbf n}^{\iota}}\restriction (n+1))$
(recalling $\underset{\sim}{\eta_{\mathbf n}^{\iota}}$ is the generic).
Clearly, $\Vdash_{\mathbb{Q}_{\mathbf n}^{\iota}} "\underset{\sim}{\nu} \in 2^{\omega}"$.
We shall prove that it's forced to be Cohen.

$(*)$ If $p\in \mathbb{Q}_{\mathbf n}^{\iota}$ and $i=1\rightarrow 2\leq nor_{\rho}(suc_p(\rho))$
for every $\rho \in T_p$, then for every $\eta \in 2^{\omega}$,
for some $\rho \in T_p$, $lg(\rho)=lg(tr(p))+m$ and if $lg(tr(p))\leq i<tr(p)+m$
then $p^{[\rho]}\Vdash "\underset{\sim}{\nu}(i)=\eta(i)"$.

We prove it by induction on $m$. For $m=1$, as $|suc_{T_{\mathbf n}}(tr(p))\setminus suc_p(tr(p))|<\frac{|suc_{T_{\mathbf n}}(tr(p))|}{2}-1$
(by clause $(g)$ of definition $2$) and for every $i\in \{0,1\}$
we have $|g_{tr(p)}^{-1}\{i\}|>\frac{\lambda_{tr(p)}}{2}-1$, hence
there are $\rho_0,\rho_1 \in suc_p(tr(p))\setminus \{\rho\}$ such
that $g_{tr(p)}(\rho_0)=0, g_{tr(p)}(\rho_1)=1$ and by the definition
of $\underset{\sim}{\nu}$, $p^{[\rho_0]}\Vdash "\underset{\sim}{\nu}(tr(p)+1)=0"$
and $p^{[\rho_1]}\Vdash "\underset{\sim}{\nu}(tr(p)+1)=1"$. Suppose
that we proved the theorem for $m$, then for some $\rho \in T_p$
of length $lg(tr(p))+m$ the conclusion holds. Now repeat the argument
of the first step of the induction for $p^{[\leq \rho]}$ to obtain
$\rho \leq \rho'$ of length $lg(tr(p))+m+1$ as required.

By $(*)$, $\underset{\sim}{\nu}$ is forced to lie in every open
dense set, hence it's Cohen.

\noindent \begin{flushright}
$\square$
\par\end{flushright}

Although the following result will not be used in the rest of the paper, it exhibits a natural property of the forcings that is of independent interest.

\textbf{Claim} \textbf{20}: If A) then B) where

A) (a) $p_i \in \mathbb{Q}_{\mathbf n}^{\iota}$ for $i<m$.

(b) $tr(p_i)=\rho$ for $i<m$.

(c) If $\iota \in \{0,1\}$ then $2\leq nor(p_i)$ for every $i<m$.

(d) If $\iota=2$ then $2\leq nor(p_i)$ for every $i<m$.

(e) $lg(\rho)<m_*<m$.

(f) There is $\rho<\eta \in T_{\mathbf n}$ such that $\lambda_{<\eta} \leq m_*<m \leq \mu_{\eta}$
(for example, it follows from the assumption $m\leq \mu_{\eta} \iff m_* \leq \lambda_{\leq \eta}$).

B) There is an equivalence relation $E$ on $\{0,1,...,m-1\}$ with
$\leq m_*$ equivalence classes such that if $i<m$ then $\{p_j : j\in (i/E)\}$
has a common upper bound.

\textbf{Proof}: Let $\eta \in T_{\mathbf n}^{[\rho \leq]}$ be as in
clause (f). Let $k_*=lg(\eta)$ and define $\lambda_{\mathbf{n},k}:=\Pi \{\lambda_{\nu} : \nu \in T_{\mathbf n}, lg(\nu)<k\}$,
$T_{\mathbf n, \rho, k}:=\{\nu \in T_{\mathbf n} : \rho \leq \nu \in T_{\mathbf n},lg(\nu)=k\}$.
Recall that $\lambda_{\nu}$ is the size of $suc_{\mathbf n}(\nu)$,
hence $|T_{\mathbf n, \rho, k_*}|$ is the product of all $\lambda_{\nu}$
such that $\rho \leq \nu$ and $lg(\nu)<k_*$, which is $\leq \lambda_{\mathbf{n},k_*}$.
For each $i<m$ let $\rho_i \in p_i$ be of length $k_*$, then $\rho_i \in T_{\mathbf n, \rho, k_*}$
by the definition of $T_{\mathbf n, \rho, k_*}$ and the assumptions
on $p_i$. Define $\rho_{i}^{+}$ for $i<m$ as follows: if $\lambda_{\eta} < \lambda_{\rho_i}$,
define $\rho_{i}^{+}:= \rho_i$. Otherwise we let $\rho_{i}^{+}\in suc_{p_i}(\rho_i)$.

Define the equivalence relation $E:=\{(i,j) : \rho_{i}^{+}=\rho_{j}^{+}\}$.
Let $j<m$, for every $i\in (j/E)$ define $p_i'=p_{i}^{[\rho_{j}^{+}]}$
(this is well defined, as $\rho_{i}^{+}=\rho_{j}^{+}$), then $tr(p_i')=\rho_{j}^{+}$
for every $i \in (j/E)$. By the choice of $\eta$, for $j<m$, $|j/E|\leq m\leq \mu_{\eta} \leq \mu_{\rho_j^+}$
(by the choice of $\rho_j^+$ and definition $2$).

By claim 10, the set $\{p_i' : i \in (j/E)\}$ has a common upper
bound, hence $\{p_i : i \in (j/E)\}$ has a common upper bound.

By the choice of $p_i^+$, the number of $E$-equivalence classes
is bounded by $\lambda_{<\eta}$. As $\lambda_{<\eta} \leq m_*$,
we're done.

\noindent \begin{flushright}
$\square$
\par\end{flushright}

\textbf{\large 6. Not adding an unwanted real}{\large \par}

A crucial step towards our final goal is to prove that the only generic
reals in finite length iterations of $\mathbb{Q}_{\mathbf n}^2$ are
the $\eta_t$s. This will be used later in order to show that $\omega^{\omega} \setminus \{\eta_t : t\in L\}$
is null with resepect to the relevant ideal. We intend to strengthen
this result dealing with arbitrary length iterations in {[}F1424{]}.

\textbf{Claim} \textbf{21}: We have $p_* \Vdash_{\mathbb{P}} "\underset{\sim}{\rho}$
is not $(\mathbb{Q}_{\mathbf n}^{\iota}, \underset{\sim}{\eta_{\mathbf n}^{\iota}})$-generic
over $V"$ when:

a) $\iota \in \{1,2\}$ and $\alpha_* <\omega$.

b) $(\mathbb{P}_{\alpha}, \underset{\sim}{\mathbb{Q}_{\alpha}} : \alpha<\alpha_*)$
is a FS iteration with limit $\mathbb{P}=\mathbb{P}_{\alpha_*}$.

c) $\mathbf{n}_{\alpha} \in \mathbf{N}$ is special (note: $\mathbf{n}_{\alpha}$
is not a $\mathbb{P}_{\alpha}-$name).

d) $\Vdash_{\mathbb{P}_{\alpha}} "\underset{\sim}{\mathbb{Q}_{\alpha}}={(\mathbb{Q}_{\mathbf{n}_{\alpha}}^{\iota})}^{V^{\mathbb{P}_{\alpha}}}"$.

e) $\mathbf{n} \in \mathbf{N}$ is special.

f) For every $\alpha$, $\mathbf{n}$ and $\mathbf{n}_{\alpha}$ are far
(i.e. $\eta_1 \in T_{\mathbf n} \wedge \eta_2 \in T_{\mathbf{n}_{\alpha}}\rightarrow \lambda_{\eta_1}^{\mathbf n} \ll \mu_{\eta_2}^{\mathbf{n}_{\alpha}}$
or $\lambda_{\eta_2}^{\mathbf{n}_{\alpha}} \ll \mu_{\eta_1}^{\mathbf n})$. Moreover, for every $\alpha<\alpha_*$ for every $l$ large enough, for
some $m\in \{l,l+1\}$ we have:

If $\rho \in T_{\mathbf n}, lg(\rho)=l$, $\nu_1 , \nu_2 \in T_{\mathbf{n}_{\alpha(l)}}$
and $lg(\nu_1) < m \leq lg(\nu_2))$ then $\lambda_{\mathbf{n}_{\alpha(l)},\nu_1} \ll \mu_{\mathbf{n}, \rho}$
and $\lambda_{\mathbf{n}, \rho} \ll \mu_{\mathbf{n}_{\alpha(l)}, \nu_2}$.

g) $p_* \Vdash_{\mathbb{P}} "\underset{\sim}{\rho} \in lim(T_{\mathbf n})"$.

\textbf{Proof}: For $\eta\in T_{\mathbf n}$ define $W_{\mathbf n,\eta}:=\{w : w\subseteq suc_{T_{\mathbf n}}(\eta)$
and $i=1\rightarrow lg(\eta) \leq nor_{\eta}^{\mathbf n}(w)$ and $i=2 \rightarrow 2\leq nor_{\eta}^{\mathbf n}(w)\}$.
For $n<\omega$ define $\Lambda_n=\{\eta \in T_{\mathbf n} : lg(\eta)<n\}$,
so $T_{\mathbf n}=\underset{n<\omega}{\cup}\Lambda_n$. Define $S_n:=\{\bar{w} : \bar{w}=(w_{\eta} : \eta \in \Lambda_n \wedge w_{\eta} \in W_{\mathbf n,\eta})\}$
and $S=\underset{n<\omega}{\cup}S_n$. $(S, \leq)$ is a tree with
$\omega$ levels such that each level is finite and $lim(S)=\{\bar{w} : \bar{w}=(w_{\eta} : \eta \in T_{\mathbf n})$
and $\bar{w} \restriction \Lambda_n \in S_n$ for every $n\}$. For
$\bar{w} \in lim(S)$ let $\mathbf{B}_{\bar{w}}:=\{ \rho \in lim(T_{\mathbf n}):$
for every $n$ large enough, $\rho \restriction (n+1) \in w_{\rho \restriction n}\}$,
so $\mathbf{B}_{\bar{w}}=\underset{m<\omega}{\cup}\mathbf{B}_{\bar{w},m}$
where $\mathbf{B}_{\bar{w},m}=\{\rho \in lim(T_{\mathbf n}) :$ if $m\leq n$
then $\rho \restriction (n+1) \in w_{\rho \restriction n}\}$. We
shall prove that 

$(*) \Vdash_{\mathbb{Q}_{\mathbf n}^{\iota}} "\underset{\sim}{\eta_{\mathbf n}^{\iota}} \in \mathbf{B}_{\bar{w}}"$
for every $\bar{w} \in lim(S)$. In fact, for every $p\in \mathbb{Q}_{\mathbf n}^{\iota}$ there is a stronger $q$ and $m<\omega$ such that $lim(q) \subseteq \bold{B}_{\bar w, m}$.

Let $p\in \mathbb{Q}_{\mathbf n}^{\iota}$, we shall prove that for
some $p\leq q$ and $m<\omega$, $q\Vdash \underset{\sim}{\eta_{\mathbf n}^{i}} \in \mathbf{B}_{\bar{w},m}$.
Let $\nu \in T_p$ such that $lg(\nu)$ is large enough and let $m=lg(\nu)$.
Now $q$ will be defined by taking the subtree obtained from the intersection
of $T_p^{[\leq \nu]}$ with $(\underset{\nu \leq \rho}{\cup}w_{\rho})$.
By the co-bigness property, $q$ is a well defined condition, and
obviously $q\Vdash \underset{\sim}{\eta_{\mathbf n}^{i}} \in \mathbf{B}_{\bar{w},m}$.

By $(*)$ it suffices to prove that for some $\bar{w} \in lim(S)$,
$p_* \nVdash_{\mathbb{P}} "\underset{\sim}{\rho} \in \mathbf{B}_{\bar{w}}"$.

Proof: Assume towards contradiction that $p\Vdash "\underset{\sim}{\rho} \in \mathbf{B}_{\bar{w}}$
for every $\bar{w} \in lim(S)"$, so there is a sequence $(p_{\bar{w}} : \bar{w} \in lim(S))$
and a sequence $(m(\bar{w}) : \bar{w} \in lim(S))$ such that:

a) $p_* \leq p_{\bar{w}}$.

b) $p_{\bar{w}} \Vdash \underset{\sim}{\rho} \in \mathbf{B}_{\bar{w},m(\bar{w})}$. 

By increasing the conditions $p_{\bar{w}}$ if necessary, we may assume
WLOG that:

1. $tr(p_{\bar{w}}(\alpha))$ is an object for every $\bar{w}$ and
every $\alpha \in Dom(p_{\bar{w}})$.

2. If $\iota=1$ and $\alpha \in Dom(p_{\bar{w}})$, then $p_{\bar w} \restriction \alpha \Vdash_{\mathbb{P}_{\alpha}} "\nu \in p_{\bar w}(\alpha) \rightarrow nor_{\nu}(Suc_{p_{\bar w}(\alpha)}(\nu))>2"$.

If $\iota=2$ and $\alpha \in Dom(p_{\bar{w}})$, then for some $m \ll lg(tr(p_{\bar{w}}(\alpha)))$,
$p_{\bar{w}} \restriction \alpha \Vdash_{\mathbb{P}_{\alpha}} \nu \in p_{\bar{w}}(\alpha) \rightarrow 1+\frac{1}{m} \leq nor(suc_{p_{\bar{w}}(\alpha)}(\nu))$.

In order to prove (1)+(2), we shall prove by induction on $\beta \leq \alpha_*$
that for every $p \in \mathbb{P}_{\beta}$ there is $p\leq q \in \mathbb{P}_{\beta}$
satisfying (2) and forcing a value to the relevant trunks. 

The induction step: assume that $\beta=\gamma+1$. As $p(\gamma)$
is a $\mathbb{P}_{\gamma}-$name of a condition in $\mathbb{Q}_{\mathbf n}^2$,
there are $p \restriction \gamma \leq p' \in \mathbb{P}_{\gamma}$
and $\rho$ such that $p' \Vdash_{\mathbb{P}_{\gamma}} tr(p(\gamma))=\rho$.
As $p' \Vdash_{\mathbb{P}_{\gamma}} p(\gamma) \in \mathbb{Q}_{\mathbf n}^2$
and by the definition of $\mathbb{Q}_{\mathbf n}^2$, there is $p' \leq p''$
and $m\leq \mu_{lg(\rho)}$ such that $p'' \Vdash_{\mathbb{P}_{\gamma}} \nu \in p(\gamma) \rightarrow 1+\frac{1}{m} \leq nor(suc_{p(\gamma)}(\nu))$.
Now choose $m \ll m_1$, so $p'' \Vdash_{\mathbb{P}_{\gamma}} "$there
is $\nu \in p(\gamma)$ such that $lg(\nu)=m_1"$. Therefore there
are $p''\leq p^*$ and $\nu$ of length $m_1$ such that $p^* \Vdash_{\mathbb{P}_{\gamma}} "\nu \in p(\gamma) \wedge (\nu \leq \eta \in p(\gamma) \rightarrow 1+\frac{1}{m} \leq nor(suc_{p(\gamma)}(\eta)))"$.
By the induction hypothesis, there is $p^* \leq q' \in \mathbb{P}_{\gamma}$
satisfying (1)+(2). Now define $q:=q' \cup (\gamma,p(\gamma)^{[\nu \leq]})$,
obviously $q$ is as required. The proof for $\mathbb{Q}_{\mathbf n}^1$
is similar.

Now we shall define a partition of $lim(S)$ to $\aleph_0$ sets as
follows:

Let $W_{m,u,\bar{\rho}}=\{\bar{w} \in lim(S) : m(\bar{w})=m, Dom(p_{\bar{w}})=u\in [\alpha_*]^{<\aleph_0}, \bar{\rho}=(tr(p_{\bar{w}}(\alpha)) : \alpha \in u)\}$.
Choose $(m_*, u_*, \bar{\rho}_*)$ such that $W=W_{m_*, u_*, \bar{\rho}_*}\subseteq lim(S)$
is not meagre. Let $\bar{u}_* \in S$ such that $W$ is comeager
above $\bar{u}_*$. Let $l(*)$ be such that $\bar{u}_* \in S_{l(*)}$.

Denote $\bar{\rho}^*=(\rho_{\alpha}^{*} : \alpha \in u_*)$, let $(\alpha_n : n<n(*))$
list $u_*$ in increasing order and let $\alpha_{n(*)}=\alpha_*$.
Therefore, if $\bar{u}_* \leq \bar{w} \in W$ then $Dom(p_{\bar{w}})=\{\alpha_0,...,\alpha_{n(*)-1}\}$
and $tr(p_{\bar{w}}(\alpha_n))=\rho_{\alpha_n}^*$ for every $n<n(*)$.

By our assumption, $\mathbf n$ is far from $\mathbf{n}_{\alpha}$. As
increasing $\bar{u}_*$ is not going to change the argument, we may
assume that $l(*)$ is large enough so $\underset{\alpha \in u_*}{\wedge}lg(\rho_{\alpha}^*)<l(*)$
and if $l<n(*)$, $\nu \in T_{\mathbf n}$, $\rho \in T_{\mathbf{n}_{\alpha_l}}$
and $lg(\bar{u}_*) \leq lg(\nu)$, then $\lambda_{\mathbf n,\nu} \ll \mu_{\mathbf{n}_{\alpha_l},\rho}$
or $\lambda_{\mathbf{n}_{\alpha_l},\rho} \ll \mu_{\mathbf n,\nu}$. Note
that we don't have to assume that $lg(\bar{u}_*)\leq lg(\rho)$: For
every $n<n(*)$, there is $m_n$ as guaranteed by $(f)$, with
$(\mathbf{n}_{\alpha_n}, lg(\nu),m_n)$ here standing for $(\mathbf{n}_{\alpha},l,m)$
there. If $lg(\rho) \leq m_n$, then by taking an arbitrary $\nu_2$
of length $>m_n$, it follows from $(f)$ that $\lambda_{\mathbf{n}_{\alpha_n},\rho} \ll \mu_{\mathbf n,\nu}$.
If $m_n< lg(\rho)$, then by taking an arbitrary $\nu_2$ of length
$\leq m_n$, we get $\lambda_{\mathbf{n},\nu} \ll \mu_{\mathbf{n}_{\alpha_n},\rho}$.

Recalling $(f)$ (and by increasing $\bar{u}_*$ if necessary),
let $(m_n : n<n(*))$ be a series of natural numbers such that $(\mathbf n, \mathbf{n}_{\alpha_n}, l(*), m_n)$
satisfy that assumptions of $(f)$ (with $(\mathbf n, \mathbf{n}_{\alpha_n}, l(*), m_n)$
here standing for $(\mathbf n,\mathbf{n}_{\alpha},l,m)$ there).

Let $\Lambda_m^0=\Lambda_{m+1} \setminus \Lambda_m=\{ \rho \in T_{\mathbf n} : lg(\rho)=m\}$
and let $S_m^0=\{\bar{w} : \bar{w}=(w_{\eta} : \eta \in \Lambda_m^0),$
for every $\eta \in \Lambda_m^0,$ $w_{\eta} \in W_{\mathbf n,\eta} \}$.

Recalling that above $\bar{u}_*$, $W$ is nowhere meagre, for every
$\bar{v} \in S_{l(*)}^0$ there is $\bar{w}_{\bar{v}} \in W \subseteq lim(\mathbb{S})$
such that $\bar{u}_*{}^\frown \bar{v} \leq \bar{w}_{\bar{v}}$.

Choose $p_n, U_n$ by induction on $n\leq n(*)$ such that:

1. $p_n \in \mathbb{P}_{\alpha_n}$.

2. If $m<n$ then $p_m \leq p_n \restriction \alpha_m$.

3. $U_n \subseteq S_{l(*)}^0$.

4. If $m<n$ then $U_n \subseteq U_m$.

5. If $E$ is an equivalence relation on $U_n$ with $\leq \Pi \{|T_{\mathbf{n}_{\alpha_l},m_l}| : n\leq l<n(*)\}$
equivalence classes, then for some $\bar{v}_* \in U_n$, $\cap \{ \underset{\rho \in T_{\mathbf n,l(*)}}{\cup} w_{\bar{v}_\rho} : \bar{v} \in \bar{v}_*/E\}=\emptyset$.

6. If $\bar{v} \in U_n$ then $p_{\bar{w}_{\bar{v}}} \restriction \alpha_n \leq p_n$.

Suppose we've carried the induction, then for every $\bar{v} \in U_{n(*)}$,
$p_{\bar{w}_{\bar{v}}}=p_{\bar{w}_{\bar{v} \restriction \alpha_{n(*)}}} \leq p_{n(*)}$,
hence by the choice of $p_{\bar{w}_{\bar{v}}}$, $p_{n(*)} \Vdash \underset{\sim}{\rho} \in \cap \{\mathbf{B}_{\bar{w}_{\bar{v}}, m_*} : \bar{v} \in U_{n(*)}\}$.
Therefore it's enough to show that $\cap \{\mathbf{B}_{\bar{w}_{\bar{v}}, m_*} : \bar{v} \in U_{n(*)}\}=\emptyset$.
By its definition, $\mathbf{B}_{\bar{w}_{\bar{v}},m_*}=lim(T_{\bar{v}})$
where $T_{\bar{v}}=\{ \eta \in T_{\mathbf n} : $ if $m_*<lg(\eta)$
then $\eta(m+1) \in w_{\eta \restriction m}$ for every $m_* \leq m \}$.
Therefore, if we show that $\cap \{T_{\bar{v}} \cap T_{\mathbf n,l(*)+1} : \bar{v} \in U_{n(*)}\}=\emptyset$,
then it will follow that $\cap \{lim(T_{\bar{v}}) : \bar{v} \in U_{n(*)}\}=\emptyset$.
This follows from part (5) of the induction hypothesis, as $\cap \{ \underset{\rho \in T_{\mathbf n,l(*)}}{\cup} w_{\bar{v}_\rho} : \bar{v} \in U_{n(*)} \}=\emptyset$.
This contradiction proves the claim.

Carrying the induction: For $n=0$, choose any $p_0 \in \mathbb{P}_{\alpha_0}$
and let $U_0=S_{l(*)}^0$. It's enough to show that $U_0$ satisfies
$(5)$. Let $E$ be an equivalence relation on $U_0$ with $m_{**} \leq \Pi \{|T_{\mathbf{n}_{\alpha(l)},m_l}| : l<n(*)\}$
equivalence classes and denote $\Pi \{|T_{\mathbf{n}_{\alpha(l)},m_l}| : l<n(*)\}$
by $m'$. For every $m<m_{**}$, denote by $U_{0,m}$ the $m$th
equivalence class of $E$. Suppose towards contradiction that for
every $m<m_{**}$ there is some $\eta_m$ in $\cap \{\underset{\rho}{\cup}w_{\rho} : \bar{w} \in U_{0,m}\}$.
For every $m$ there is $\rho_m$ such that $\eta_m \in suc_{T_{\mathbf n}}(\rho_m)$.
Choose $\bar{w}=(w_{\rho} : \rho \in T_{\mathbf{n},l(*)})$ by letting
$w_{\rho}=suc_{T_{\mathbf n}}(\rho) \setminus \{ \eta_m : m<m_{**} \wedge \rho_m=\rho\}$.
We shall prove that $\bar{w} \in U_0$. It will then follow that $\bar{w} \in U_{0,m}$
for some $m$, therefore $\eta_m \in \underset{\rho}{\cup}w_{\rho}$,
contradicting the definition of $w_{\rho}$. This proves that $U_0$
is as required. In order to prove that $\bar{w} \in U_0$, note
that for every $\rho$, $|suc_{T_{\mathbf n}}(\rho) \setminus w_{\rho}| \leq |\{m : \rho_m=\rho\}| \leq m_{**} \leq m'=\Pi \{|T_{\mathbf{n}_{\alpha(l)},m_l}| : l<n(*)\} \ll \mu_{\mathbf n, \rho}$
(the last inequality follows by $(f)$ and the choice of $m_l$, recalling that the $m_l$ were chosen to satisfy the assumptions of (f) and recalling the definition of the $\lambda_{\bold n, \nu}$).
Therefore, $\bar{w} \in U_0$.

Suppose now that $n=k+1\leq n(*)$. Choose $q_k \in \mathbb{P}_{\alpha_k}$
such that $p_k \leq q_k$ and $q_k$ forces a value $\Lambda_{\bar{v}}^{k}$
to $\{\rho \in p_{\bar{w}_{\bar{v}}}(\alpha_k) : lg(\rho)=m_k\}$
for every $\bar{v} \in U_k$. For every $\rho \in T_{\mathbf{n}_{\alpha_k},m_k}$
let $U_{k,\rho}=\{\bar{v} \in U_k : \rho \in \Lambda_{\bar{v}}^{k}\}$.
If $\bar{v} \in U_k$, then $q_k$ forces the value $\Lambda_{\bar{v}}^{k}$
to $\{\rho \in p_{\bar{w}_{\bar{v}}}(\alpha_k) : lg(\rho)=m_k\}$,
hence $U_k=\cup \{U_{k,\rho} : \rho \in T_{\mathbf{n}_{\alpha},m_k}\}$.
WLOG $U_{k,\rho}$ are pairwise disjoint. Now suppose towards contradiction
that none of them satisfies requirement (5) of the induction for $k+1$,
then each $U_{k,\rho}$ has a counterexample $E_{\rho}$, and the
union $\underset{\rho}{\cup}E_{\rho}$ is therefore an equivalence
relation which is a counterexample to $U_k$ satisfying (5). Therefore,
for some $\rho$, $U_{k,\rho}$ satisfies (5), so choose $U_n=U_{k,\rho}$.

Define $p_n \in \mathbb{P}_{\alpha_k +1} \subseteq \mathbb{P}_{\alpha_n}$
as follows:

1. $p_n \restriction \alpha_k =q_k$.

2. $p_n(\alpha_k)=\cap \{p_{\bar{w}_{\bar{v}}}(\alpha_k)^{[\rho \leq]} : \bar{v} \in U_n\}$.

Now for every $\bar{v} \in U_k$, $p_{\bar{w}_{\bar{v}}} \restriction \alpha_k \leq p_k \leq q_k$,
hence $q_k \Vdash_{\mathbb{P}_{\alpha_k}} \nu \in p_{\bar{w}_{\bar{v}}}(\alpha_k) \rightarrow 1+\frac{1}{m} \leq nor(suc_{p_{\bar{w}_{\bar{v}}}(\alpha_k)}(\nu))$.
We shall prove that $q_k \Vdash_{\mathbb{P}_{\alpha_k}} p_n(\alpha_k) \in \mathbb{Q}_{\mathbf{n}_{\alpha}}^2$.
As, $|U_n| \leq |S_{l(*)}^0| \leq 2^{\Sigma\{\lambda_{\mathbf n,\rho'} : \rho' \in \Lambda_{l(*)}^0\}}<\mu_{\mathbf{n}_{\alpha_k},\rho}$ (with the last inequality following from (f) and the choice of the $m_k$'s),
the assumptions of claim 10 hold, the conclusion follows by the proof
of claim 10. A similar argument (using the first part of claim 10)
proves the claim for the case of $\mathbb{Q}_{\mathbf n}^1$.

So $p_n$ obviously satisfies requirements 1,2 and 6.

\noindent \begin{flushright}
$\square$
\par\end{flushright}

\textbf{\large 7. Main measurability claim}{\large \par}

We're now ready to prove the main result. We shall first prove that
Cohen forcing (hence $\mathbb{Q}_{\mathbf n}^i$) turns the ground model
set of reals into a null set with respect to our ideal. We will then
prove the main result by using a Solovay-type argument.

\textbf{Claim} \textbf{22}: For $\iota \in \{0,\frac{1}{2},1,2\}$
we have $\Vdash_{Cohen} "$there is a Borel set $\mathbf{B} \subseteq lim(T_{\mathbf{n}_*})$
such that $lim(T_{\mathbf{n}_*})^V \subseteq \mathbf{B}$ and $\mathbf{B}$
is $(\mathbb{Q}_{\mathbf{n}_*}^{\iota}, \underset{\sim}{\eta_{\mathbf{n}_*}^{\iota}})$-null$"$ (where by "$(\mathbb{Q}_{\mathbf{n}_*}^{\iota}, \underset{\sim}{\eta_{\mathbf{n}_*}^{\iota}})$-null" we mean that for every $p$ there is a stronger $q$ with $lim(q) \cap \bold{B}=\emptyset$)

\textbf{Proof}: Let $\mathbb{Q}$ be the set of finite functions with
domain $\{\eta\in T_{\mathbf{n}_*} : lg(\eta)<k\}$ for some $k<\omega$
such that $f(\rho) \in suc_{T_{\mathbf{n}_*}}(\rho)$. $(\mathbb{Q}, \subseteq)$
is countable and for every $q\in \mathbb Q$ there are $q\leq q_1,q_2 \in \mathbb Q$
which are incompatibe, hence is equivalent to Cohen forcing. Let $\underset{\sim}{f}:=\underset{g\in \underset{\sim}{G}}{\cup}g$.
For $f\in S=\Pi \{suc_{T_{\mathbf{n}_*}}(\rho) : \rho \in T_{\mathbf{n}_*}\}$
define $\mathbf{B}_f:= \{\eta \in lim(T_{\mathbf{n}_*}) :$ for infinitely
many $n$ we have $\eta \restriction (n+1)=f(\eta_n)\}$. For every
$n<\omega$ let $\mathbf{B}_{f,n}=\{\eta \in lim(T_{\mathbf{n}_*}) : \eta \restriction (m+1)\neq f(\rho)$
if $n\leq m$ and $n\leq lg(\rho)\}$. Clearly, $\Vdash "\underset{\sim}{f} \in S"$,
$\mathbf{B}_f^c=\underset{n<\omega}{\cup}\mathbf{B}_{f,n}$, and obviously
each $\mathbf{B}_{f,n}$ is Borel, hence $\mathbf{B}_f$ is Borel. For
every $\eta \in T_{\mathbf{n}_*}$ let $w_{\eta}=suc_{T_{\mathbf{n}_*}}(\eta)\setminus \{f(\eta)\}$.
As in claim 21, $\Vdash_{\mathbb{Q}_{\mathbf{n}_*}^{\iota}} "\underset{\sim}{\eta_{\mathbf{n}_*}^{\iota}} \in \mathbf{B}_{\bar{w}}"$
for $\bar{w}$ and $\mathbf{B}_{\bar{w}}$ as in that proof. Hence $\Vdash_{\mathbb{Q}_{\mathbf{n}_*}^{\iota}} "\underset{\sim}{\eta_{\mathbf{n}_*}^{\iota}} \notin \mathbf{B}_f$,
so $\mathbf{B}_f$ is $(\mathbb{Q}_{\mathbf{n}_*}^{\iota}, \underset{\sim}{\eta_{\mathbf{n}_*}^{\iota}})$-null.
Let $G\subseteq \mathbb Q$ be generic and let $g=\underset{\sim}{f}[G]$,
so $\mathbf{B}_g$ is a $(\mathbb{Q}_{\mathbf{n}_*}^{\iota}, \underset{\sim}{\eta_{\mathbf{n}_*}^{\iota}})$-null
Borel set in $V[G]$. We shall prove that $V[G] \models lim(T_{\mathbf{n}_*})^V \subseteq \mathbf{B}_g$.
Let $\eta \in lim(T_{\mathbf{n}_*})^V$ and $m<\omega$, it's enough
to show that in $V$, $\Vdash_{\mathbb Q} "$for some $m\leq k$ and
$\rho \in T_{\mathbf{n}_*}, \underset{\sim}{f}(\rho)=\eta \restriction(k+1)"$.
Let $p\in \mathbb{Q}$, we can extend $p$ to a function $p\leq q$
with domain $\{\eta \in T_{\mathbf{n}_*} : lg(\eta)<k\}$ for some $m\leq k$.
Now let $q\leq s$ be an extension of $q$ with domain $\{\eta \in T_{\mathbf{n}_*} : lg(\eta) \leq k\}$
such that $s(\eta \restriction k)=\eta \restriction (k+1)$. Obviously,
$s$ forces the required conclusion, so we're done.

\noindent \begin{flushright}
$\square$
\par\end{flushright}

\textbf{Main} \textbf{conclusion} \textbf{23}: Let $i\in \{1,2\}$.
Let $V\models CH$ and suppose $\aleph_1 < \mu=\mu^{\aleph_0}$.
Let $L$ be a linear order of cardinality $\mu$ that is homogeneous, i.e. that any two nonempty open intervals are isomorphic (for an example of such a linear order, see e.g. Section 4 in [KrSh859])  \footnote{Note that such an order $L$ is dense with no endpoints, and that if $-\infty=s_0^l<s_1^l<...<s_{n-1}^l<s_n^l=\infty$  ($l=0,1$), then there is an automorphism $\pi$ of $L$ such that $\pi(s_k^0)=s_k^1$. In addition, if $s_k^0=s_k^1$ and $s_{k+1}^0=s_{k+1}^1$, then $\pi$ can be the identity on $(s_k^0, s_{k+1}^0)$.}.
Suppose that $\mathbf q$ is as in 13(A) such that $L_{\mathbf q}=L$
and $\mathbf{m}_t=\mathbf{m}$ for every $t\in L_{\mathbf q}$ is a (constant)
definition of the forcing $\mathbb{Q}_{\mathbf n}^i$, then:

a) $\mathbb{P}_{\mathbf q}$ is a c.c.c. forcing notion of cardinality
$\mu$.

b) $\Vdash_{\mathbb{P}_{\mathbf q}} "2^{\aleph_0}=\mu"$.

c) Let $G\subseteq \mathbb{P}_{\mathbf q}$ be generic over $V$ and let $\eta_t=\underset{\sim}{\eta_t}[G]$
for $t\in L_{\mathbf q}$. In $V[G]$ we have the sets $X:=\{\eta_t : t\in L_{\mathbf q}\}$ and $<_X:=\{(\eta_s, \eta_t) : s<_{L_{\bold q}} t\}$. Note that these sets are definable over $V$: $\eta \in X$ iff it satisfies "$\eta$ is $(\mathbb{Q}_{\mathbf n}^i, \underset{\sim}{\eta})$-generic over $V$", and $(\eta_s, \eta_t) \in <_X$ iff they satisfy "$\eta_s$ is not $(\mathbb{Q}_{\mathbf n}^i, \underset{\sim}{\eta})$-generic over $V[\eta_t]$". Now let $V[X^+]$ be the collection of sets hereditarily definable from elements of $V$ and finite sequences of members of $X^+:=X\cup \{X, <_X\}$, so $X, <_X \in V[X^+]$ \footnote{The addition of ${X, <_X}$ was done only for the sake of clarity. We could have worked instead in $V[X]$, i.e. the collection of sets hereditarily definable from finite sequences of members of $X$}. Similarly, for $Z\subseteq X$, let $Z^+:=Z\cup \{X,<_X\}$. Note that, in $V[G]$, if $y\subseteq H(\aleph_0)$ then $y\in V[X^+]$ iff $y\in V[Z^+]$ for some finite $Z\subseteq X$.

$(\alpha)$ $V[X^+]\models ZF+\neg AC_{\aleph_0}$ and $lim(T_{\mathbf{n}})^{V[X^+]}=\cup \{lim(T_{\mathbf{n}})^{V[\{\eta_t : t \in u\}]} : u\subseteq L_{\mathbf q}$
is finite$\}$.

$(\beta)$ $(\mathbb{Q}_{\mathbf n}^i,\aleph_1)$-measurability holds in $V[X^+]$: Every
$A\subseteq lim(T_{\mathbf n})^{V[X^+]}$ is $I_{\mathbb{Q}_{\mathbf n}^i,\aleph_1}$-measurable.

$(\gamma)$ $V[X^+] \models "\{\eta_t : t\in L_{\mathbf q}\}=lim(T_{\mathbf{n}})$ $mod$
$I_{\mathbb{Q}_{\mathbf n}^i,\aleph_1}"$.

$(\delta)$ In $V[X^+]$, if $J\subseteq L_{\mathbf q}$ is a proper initial segment
then $\{\eta_t : t\in J\}\in I_{\mathbb{Q}_{\mathbf n}^i,\aleph_1}$.

$(\epsilon)$ In $V[X^+]$, the ideal $I_{\mathbb{Q}_{\mathbf n}^i,\aleph_1}$ is
non-trivial.

$(\zeta)$ $\aleph_1$ is not collapsed, there is an $\omega_1$-sequence
of different reals, and if $V=L$ (here $L$ is the constructible universe) then $\aleph_1^{L}=\aleph_1^{V[X^+]}$.

\textbf{Proof}: Clause a) By the definition of $\mathbb{P}_{\mathbf q}$
and claim 17, so $|\mathbb{P}_{\mathbf q}| \leq \Sigma\{ |\mathbb{P}_{\mathbf q,J}| : J\subseteq L$
is finite$\} \leq 2^{\aleph_0}+|L|^{<\aleph_0}=2^{\aleph_0}+\mu=\mu$.
\\
\\
Clause b) By a) we have $\Vdash_{\mathbb{P}_{\mathbf q}} "2^{\aleph_0} \leq \mu"$,
and as $|L|=\mu$ we have $\Vdash_{\mathbb{P}_{\mathbf q}} "\mu=|L|\leq |\{\underset{\sim}{\eta_t} : t\in L\}| \leq 2^{\aleph_0}"$.
Together we're done.
\\
\\
Clause c) $(\alpha)$ By the definitions of $V[X^+]$ and $\mathbb{P}_{\mathbf q}$.
In particular, $\neg AC_{\aleph_0}$, as we can use $(A_n : n<\omega)$
where $A_n:=\{\{\eta_{t_l} : l<n\} : t_0<_L...<_L t_{n-1}\}$. As $V[X^+]$ is really just $HOD(V\cup X^{<\omega})$ in $V[G]$, $V[X^+] \models ZF$ follows by the standard arguments in the literature.
\\
\\
Clause c)$(\beta)$ Let $A\in V[X^+]$ be a subset of $lim(T_{\mathbf{m}_*})$.
$A$ is definable in $V[G]$ by a first order formula $\phi(x, \bar{a},{c})$
such that $c\in V$ and $\bar{a}=(\eta_{t_0},...,\eta_{t_{n-1}})$
is a finite sequence from $X$. Let $J=\{s\in L_{\mathbf q} : s\leq t_l$
for some $l\}$. For $s\in L\setminus J$ let $L_s=\{t_l : l<n\} \cup \{s\}$,
then $L_s \in \mathbf{L}_{\mathbf q}$ hence by 14 we have $\mathbb{P}_{L_s}\lessdot \mathbb{P}_{L_{\mathbf q}}$.
Let $\underset{\sim}{T_s}=TV(\phi(\underset{\sim}{\eta_s},\bar{a},c))$,
so $\underset{\sim}{T_s}$ is a $\mathbb{P}_{L_{\mathbf q}}$-name and
actually a $\mathbb{P}_{L_s}$-name.

Let $(p_{s,i} : i<\omega)$ be a maximal antichian in $\mathbb{P}_{L_s}$
and let $W_s \subseteq \omega$ such that $p_{s,i}\Vdash \underset{\sim}{T_s}=true$
if and only if $i\in W_s$. Define the $\mathbb{P}_{\{t_l : l<n\}}$-name
$\underset{\sim}{U}:= \{i<\omega : p_{s,i} \restriction \{t_l : l<n\} \in \underset{\sim}{G_{\mathbb{P}_{\{t_l : l<n\}}}}\}$.

If $G_0 \subseteq \mathbb{P}_{\{t_l : l<n\}}$ is generic over $V$
and $U=\underset{\sim}{U}[G_0]$, then in $V[G_0]$, $(lim(p_{s,i}(s)[G_0]) : i \in U)$
are pairwise disjoint: by claim $7$, if $p,q\in \mathbb{Q}_{\mathbf n}^{\iota}$
are incompatible and $\eta \in lim(p)$, then $\eta \notin lim(q)$
(otherwise, WLOG $lg(tr(p))\leq lg(tr(q))$, and both $tr(p)$ and
$tr(q)$ are initial segments of $\eta$, hence $tr(p)\leq tr(q) \in T_p$
which is a contradiction by claim $7$). Hence it's enough to show
that $((p_{s,i}(s)[G_0]) : i \in U)$ is an antichain in $V[G_0]$.
Assume towards contradiction that for some $i \neq j\in U$ there
is a common upper bound $q$ for $p_{s,i}(s)[G_0]$ and $p_{s,j}(s)[G_0]$.
Therefore there is a $\mathbb{P}_{\{t_l : l<n\}}$-name $\underset{\sim}{q}$
and $r\in G_0$ such that $r\Vdash_{\mathbb{P}_{\{t_l : l<n\}}} "p_{s,i}(s), p_{s,j}(s) \leq \underset{\sim}{q}"$.
Since $i,j\in U$, we have $p_{s,i}\restriction \{t_l : l<n\}, p_{s,j}\restriction \{t_l : l<n\} \in G_0$,
and as $G_0$ is directed, there is a common upper bound $r_1 \in G_0$
for $p_{s,i}\restriction \{t_l : l<n\}, p_{s,j}\restriction \{t_l : l<n\}$
and $r$. Now let $r^+:=r_1 \cup \{(s,\underset{\sim}{q})\} \in \mathbb{P}_{L_s}$,
then obviousy $r^+$ is a common upper bound (in $\mathbb{P}_{L_s}$)
for $p_{s,i}$ and $p_{s,j}$, which contradicts our assumption.

Moreover, $(p_{s,i}(s)[G_0] : i\in U)$ is a maximal antichain: If
$q\in {\mathbb{Q}_{\mathbf n}^{\iota}}^{V[G_0]}$ is incompatible with
$p_{s,i}(s)[G_0]$ for every $i\in U$, then as before, there are
$r\in G_0$ and a $\mathbb{P}_{\{t_l : l<n\}}$-name $\underset{\sim}{q}$
such that $r$ forces that $\underset{\sim}{q}$ is incompatible with
$p_{s,i}(s)$ for every $i\in U$. As before we can get a member of
$\mathbb{P}_{L_s}$ that is incompatible with $(p_{s,i} : i<\omega)$,
contradicting its maximality. Hence $(p_{s,i}(s)[G_0] : i\in U)$
is a maximal antichain in $V[G_0]$.

If $s_1,s_2 \in L_{\mathbf q}\setminus J$, by the homogeneity assumption,
there is an autommorphism $f$ of $L_{\mathbf q}$ over $J$ such that
$f(s_1)=s_2$. Therefore the natural map induced by $f$ is mapping
$\bar{a}$ to itself and $\underset{\sim}{\eta_{s_1}}$ to $\underset{\sim}{\eta_{s_2}}$.
Hence $\underset{\sim}{T_{s_1}}$ is mapped to $\underset{\sim}{T_{s_2}}$.
As $(\hat{f}(p_{s_1,i}) : i<\omega)$ and $W_{s_1}$ have the same
properties (with respect to $\underset{\sim}{T_{s_2}}$) as $(p_{s_2,i} : i<\omega)$
and $W_{s_2}$, we may assume WLOG that $W_{s_1}=W_{s_2}$ (denote
it by $W$) and $\hat{f}(p_{s_1,i})=p_{s_2,i}$.

Therefore, if $G_0 \subseteq \mathbb{P}_{\{t_l : l<n\}}$ is generic
and $i\in \underset{\sim}{U}[G_0]$, then there is $p_i \in {(\mathbb{Q}_{\mathbf n}^{\iota}})^{V[G_0]}$
and $W$ such that for every $s\in L\setminus J$, $p_{s,i}(s)[G_0]=p_i$
and $W_s=W$ (and $W$ can be found in the ground model).

Work now in $V[G_0]$: Let $B:=\cup \{lim(p_i) : i\in W \cap U\}$,
so $B$ is a Borel set and we shall prove that $A=B$ modulo the ideal:
by clauses $(c)(\gamma)+(c)(\delta)$ proved below, it's enough to
show that if $s\in L_{\mathbf q}\setminus J$, then ${\eta_s}\notin A\Delta B$ (note that, by its definition, $J\in V$ and hence $J\in V[X^+]$).

Let $s\in L_{\mathbf q}\setminus J$ and $i\in U$, then $p_{s,i} \in \mathbb{P}_{L_s}/G_0$
and by the choice of $p_{s,i}$, $p_{s,i} \Vdash_{\mathbb{P}_{L_s}/G_0} "\phi(\underset{\sim}{\eta_s}, \bar{a},c) $
iff $\underset{\sim}{T_s}=true$ iff $i\in W"$. In other words, in
$V[G_0]$ we have: $p_i \Vdash_{\mathbb{Q}_{\mathbf n}^{\iota}} "\phi(\underset{\sim}{\eta_s}, \bar{a},c)$
iff $i\in W"$. Since $(p_i : i\in U)$ is a maximal antichain, every
$G\subseteq \mathbb{Q}_{\mathbf n}^{2}$ generic over $V[G_0]$ must
contain exactly one of the $p_i$, hence in $V[G_0]:$ $\Vdash_{\mathbb{Q}_{\mathbf n}^{\iota}} "\phi(\underset{\sim}{\eta_s}, \bar{a},c)$
iff $i\in W$ for the $p_i$ such that $p_i \in \underset{\sim}{G}"$.
Now $p_{s,i}(s)=p_i \in \underset{\sim}{G}$ iff $\underset{\sim}{\eta_s} \in lim(T_{p_{s,i}(s)})=lim(T_{p_i})$,
hence we got $\Vdash_{\mathbb{Q}_{\mathbf n}^{2}} "\phi(\underset{\sim}{\eta_s}, \bar{a},c)$
iff $i\in W$ where $i$ is such that $\underset{\sim}{\eta_s} \in lim(T_{p_i})"$.
Therefore $\Vdash_{\mathbb{Q}_{\mathbf n}^{\iota}} "\underset{\sim}{\eta_s} \in A$
iff $\underset{\sim}{\eta_s}\in B"$.
\\
\\
Clause c)$(\gamma)$ If $\rho \in lim(T_{\mathbf{n}})^{V[X^+]}\setminus \{\eta_t : t\in L_{\mathbf q}\}$,
then $\rho \in lim(T_{\mathbf{n}})^{V[\{\eta_t : t\in u\} \cup  \{X, <_X\}]}$ for
some finite $u$. By claim 21, $\rho$ is not $(\mathbb{Q}_{\mathbf{n}}^{\iota}, \underset{\sim}{\eta_{\mathbf{n}}})$-generic
over $V$. Therefore, by the definition of $I_{\mathbb{Q}_{\mathbf n, \aleph_1}^i}$,
$\Vdash_{\mathbb{P}_{\mathbf q}} "lim(T_{\mathbf{n}})\setminus \{\underset{\sim}{\eta_t} : t\in L_{\mathbf q}\}\in I_{\mathbb{Q}_{\mathbf n}^i, \aleph_1}"$. This is due to the fact that being $(\mathbb{Q}_{\bold n}^i, \underset{\sim}{\eta})$-generic over $V$ means avoiding every Borel $(\mathbb{Q}_{\bold n}^i, \underset{\sim}{\eta})$-null set from $V$, and as $V\models CH$, there are $\aleph_1$-many such sets.
Why can we use claim 21? Assume that in claim 21 $\alpha_*$ is finite,
assumptions $(a)-(e)$ and $(g)$ hold and $(f)$ is replaced by $(h)$
where:

$(h)$ $p_* \Vdash_{\mathbb{P}} \underset{\sim}{\rho} \notin \{\underset{\sim}{\eta_{\alpha}} : \alpha<\alpha_*\}$.

There is a condition $p_* \leq p_{**}$ and a natural number $k$
such that $p_{**}\Vdash_{\mathbb{P}} \underset{\sim}{\rho}\restriction k \notin \{\underset{\sim}{\eta_{\alpha}}\restriction k : \alpha<\alpha_*\}$
and $p_{**}$ forces values to $\underset{\sim}{\rho}\restriction k$
and $\underset{\sim}{\eta_{\alpha}}\restriction k$ $(\alpha<\alpha_*)$,
which will be deonted by $\rho_*$ and $\eta_{\alpha}^*$ $(\alpha<\alpha_*)$.
In addition, we shall choose $k$ to be sufficiently large.

For $\mathbf n \in \mathbf N$ and $\eta \in T_{\mathbf n}$, let $\mathbf{n}^{[\eta \leq]}$
be the natural restriction of $\mathbf n$ to $T_{\mathbf n}^{[\eta \leq]}$.
Now let $\mathbf{n}_*=\mathbf{n}^{[\rho_* \leq]}$ and $\mathbf{n}_{\alpha}^{*}=\mathbf{n}_{\alpha}^{[\eta_{\alpha}^{*} \leq]}$.
By the choice of $k$, $\mathbf{n}_*$ and $\mathbf{n}_{\eta_{\alpha}}^{*}$
are far, moreover, they satisfy assumption f of claim 21, and by
iterating $\mathbb{Q}_{\mathbf{n}_{\alpha}^{*}}^{i}$ instead, we get
the desired conclusion.
\\
\\
Clause c)$(\delta)$ By claim 19, each $\mathbb{Q}_{\mathbf{m}_t}$
adds a Cohen real, hence the set of previous generics is included
in a null Borel set by claim 22. More precisely: For $t\in L$, let $V_{2, t}:=V[G\cap \mathbb{P}_{L_{<t}}]$ and let $V_{1,t}$ be the class of elements from $V_{2,t}$ that are hereditarily definable in $V[G]$ from elements of $V$, finite sequences form $\{ \eta_s : s<t\}$ and $\{ \eta_s : s<t \}$. As $\{ \eta_s : s<t\} \subseteq lim(T_{\bold n})^{V_{1,t}}$, it suffices to show that the Cohen real $\underset{\sim}{\nu_t}$ added by $\underset{\sim}{\mathbb{Q}_t}$ is Cohen over $V_{1,t}$. As $\underset{\sim}{\nu_t}$ is Cohen over $V^{\mathbb{P}_J}$ for every finitie $J\subseteq L_{<t}$, it suffices to show that every nowhere dense tree $T\in V_{1,t}$ belongs to $V^{\mathbb{P}_J}$ for some finite $J\subseteq L_{<t}$ (and so $\nu_t \in lim(T)$). Suppose then that $A=\underset{\sim}{A}[G] \in V_{1,t}$ is a real, so there are $t_0<t_1<...<t_{n-1}<t$, a formula $\phi$ and some $a\in V$ such that $A$ is definable in $V[G]$ using $\phi(x, a, \eta_{t_0},...,\eta_{t_{n-1}})$. We claim that $A\in V^{\mathbb{P}_J}$ for some finite  $J\subseteq  L_{<t}$. As $A\in V[G\cap \mathbb{P}_{L_{<t}}]$ and $\mathbb{P}_{L_{<t}} \lessdot \mathbb{P}_{\bold q}$, $\underset{\sim}{A}$ is a $\mathbb{P}_{L_{<t}}$-name and $A=\underset{\sim}{A}[G_t]$, where $G_t=G\cap \mathbb{P}_{L_{<t}}$.  
Let $p\in G_t$ force the above facts, and WLOG $\{t_0,...,t_{n-1}\} \subseteq Dom(p)$. Now if $q\in \mathbb{P}_{L_{<t}}$ is above $p$ and forces "$i\in \underset{\sim}{A}$", then this is also forced by $q\restriction Dom(p)$. In order to prove this fact, note that we can find for every $n<\omega$ an automorphism $\pi_n$ of $L_{<t}$ such that $\pi_n$ is the identity over $Dom(p)$ and such that the sets $\pi_n(Dom(q) \setminus Dom(p))$ are pairwise disjoint. Letting $q_n:=\pi_n(q)$, each $q_n$ forces "$i\in \underset{\sim}{A}$". It follows that this is forced by $q\restriction Dom(p)$ as well: Suppose not, then there is some $q'$ above $q\restriction Dom(p)$ forcing "$i \notin \underset{\sim}{A}$". But then there is some $n<\omega$ such that $(Dom(\pi_n(q)) \setminus Dom(p)) \cap Dom(q')=\emptyset$. It follows that $q'$ and $\pi_n(q)$ are compatible, a contradiction. Therefore, $q\restriction Dom(p)$ forces "$i\in \underset{\sim}{A}$". The proof for the case of "$i\notin \underset{\sim}{A}$" is similar. It follows that $\underset{\sim}{A} \in V^{\mathbb{P}_{Dom(p)}}$, as required.
\\
\\
Clause c)$(\epsilon)$ We shall prove that, in $V[X^+]$, $X=\{ \eta_t : t\in L\} \notin I_{\mathbb{Q}_{\bold n}^i, \aleph_1}$ and so the ideal is non-trivial. So let $\bar Z =(Z_{\alpha} : \alpha<\omega_1) \in V[X^+]$ be a sequence of $(\mathbb{Q}_{\bold n}^i, \underset{\sim}{\eta})$-null sets and let $Z=\underset{\alpha<\omega_1}{\cup}Z_{\alpha}$, it suffices to show that $\eta_t \notin Z$ for every large enough $t\in L$. Let $p_* \in G$ force the above-mentioned facts about $\underset{\sim}{\bar Z}$ and let $(t_l : l<n) \in L^n$ be an increasing sequence containing $Dom(p_*)$ and all $t\in L$ relevant for $\underset{\sim}{\bar Z}$. Let $t\in L$ such that $t_{n-1}<t$, we shall prove that $\eta_t \notin Z$. Let $\alpha<\omega_1$ and suppose that $\underset{\sim}{q_1}$ is a $\mathbb{P}_{\{s_l : l<k\}}$-name for a member of $(\mathbb{Q}_{\mathbf n}^i)^{V[\eta_{s_l} : l<k]}$ with $(s_l : l<k) \in (L_{<t})^k$. Then there are $p' \in \mathbb{P}_{\mathbf q}$ above $p_*$ and $\underset{\sim}{q_2}'$ such that $\underset{\sim}{q_2}'$ is a $\mathbb{P}_{\{s_l' : l<m\}}$-name, $p'$ forces "$\underset{\sim}{q_1} \leq \underset{\sim}{q_2'}$ and $lim(q_2') \cap \underset{\sim}{Z_{\alpha}}=\emptyset$" and WLOG $\{t_0,...,t_{n-1}\} \subseteq \{s_0,...,s_{k-1}\} \subseteq \{s_0',...,s_{m-1}'\}$. There is an automorphism $\pi$ of $L$ that is the identity over $\{s_0,...,s_{k-1}\}$ such that $\pi(\{s_0',...,s_{m-1}'\}) \subseteq L_{<t}$, so we may assume WLOG that $\{s_0',...,s_{m-1}'\} \subseteq L_{<t}$. It follows that $\eta_t \notin Z_{\alpha}$, and therefore, $\eta_t \notin Z$.
\\
\\
Clause c)$(\zeta)$ $V\models AC$, therefore there is an $\omega_1$-sequence
of distinct reals in $V$. $\mathbb{P}_{\mathbf q} \models ccc$, therefore
$\aleph_1$ is not collapsed, and that sequence is as required in
$V[X^+]$ as well. If $V=L$, then $\aleph_1^L=\aleph_1^{V[X^+]}$ follows
from ccc.

\noindent \begin{flushright}
$\square$
\par\end{flushright}

\textbf{\large 8. An Application to $\Pi^1_n$ Singletons}{\large \par}

We conclude the paper with an easy application of $\mathbb{Q}^2_{\mathbf n}$ that is of independent interest. By a classical result of Jensen ([Je]), there exists a forcing $\mathbb P \in L$ that adds a $\Pi^1_2$ singleton over $L$. Jensen's construction relies heavily on structural properties of $L$ such as diamond. Thanks to the explicit definability of $\mathbb{Q}^2_{\mathbf n}$ and its property of adding a unique generic real, we are able to get a $\Pi^1_2$ singleton over $L$ almost "for free". As we saw, the existence and the relevant properties of $\mathbb{Q}_{\mathbf n}^2$ are already established in $ZFC$, and the only extra assumption needed for our new construction of a $\Pi^1_2$ singleton is that the ground model reals are constructible. Our construction easily generalizes to other models of set theory, in which case if the ground model reals are  $\Sigma^1_n$ then the new singleton will be $\Pi^1_n$. Below we shall only deal with lightface definitions, so $\Sigma^1_n$ will always mean "lightface $\Sigma^1_n$". 
\\
\\
Throughout the rest of this section, fix a computable $\mathbf n \in \mathbf N$ (e.g. the one from Observation 4A).
\\
\\
\textbf{Claim 24}: Let $\mathbb Q=\mathbb{Q}^2_{\mathbf n}$ and let $G\subseteq \mathbb Q$ be generic over $V$. Suppose that $(\omega^{\omega})^V$ is lightface $\Sigma^1_n$ definable in $V[G]$, then letting $\eta$ be the canonical generic real added by $\mathbb Q$, the singleton $\{ \eta\}$ is a lightface $\Pi^1_n$ singleton.
\\
\\
\textbf{Proof}: $\eta$ is $\mathbb Q$-generic over $V$ iff 
\\
(*) For every maximal antichain $I\subseteq \mathbb Q$ from $V$, there exists $p\in I$ such that $\eta \in lim(p)$. 
\\
As "$I$ is a maximal antichain in $\mathbb Q$" is a lightface Borel statement by Claim 11, and as $(\omega^{\omega})^V$ is lightface $\Sigma^1_n$ in $V[G]$, it follows that (*) is a lightface $\Pi^1_n$ statement in $V[G]$. By the uniqueness of the generic real (Claim 21), (*) defines a $\Pi^1_n$ singleton in $V[G]$. $\square$
\\
\\
The assumptions of Claim 24 hold in $L$ for $n=2$. By the following result of Steel, canonical inner models for Woodin cardinals satisfy the assumptions for $n>2$:
\\
\\
\textbf{Theorem 25}: Let $\mathbb Q$ be a Borel ccc forcing. Let $n>2$ and let $\mathcal{M}_{n-2}$ be the least inner model with $n-2$ Woodin cardinals, then $(\omega^{\omega})^{\mathcal{M}_{n-2}}$ is lightface $\Sigma^1_n$-definable in $\mathcal{M}_{n-2}^{\mathbb Q}$.
\\
\\
\textbf{Proof}: See the proof of Theorem 3.4 in [St]. $\square$
\\
\\
By Shoenfield's absoluteness theorem, the minimal possible complexity of a nonconstructible singleton is $\Pi^1_2$. In order to obtain a similar optimality result over $\mathcal{M}_{n-2}$, we shall use the following absoluteness result due to Woodin:
\\
\\
\textbf{Theorem 26}: Let $\mathbb Q$ be a Borel ccc forcing. Let $n>2$ and let $\mathcal{M}_{n-2}$ be as in the previous theorem, then for every $\Sigma^1_n$ formula $\phi(x)$ and $a\in (\omega^{\omega})^{\mathcal{M}_{n-2}}$, $\mathcal{M}_{n-2} \models \phi(a)$ iff $\mathcal{M}_{n-2}^{\mathbb Q} \models \phi(a)$.
\\
\\
\textbf{Proof}: See e.g. Section 4 in [St] or Lemma 1.17 in [MSW]. $\square$
\\
\\
Putting everything together we get the main result of this section:
\\
\\
\textbf{Theorem 27}: Let $\mathbb Q=\mathbb{Q}_{\mathbf n}^2$.
\\
a. Suppose that $V\models ZFC+$"all reals are constructible", then $\mathbb Q$ adds a $\Pi^1_2$ singleton over $V$. In particular, $\mathbb Q$ adds a new $\Pi^1_2$ singleton over $L$ and over any forcing extension of $L$ not adding new reals.
\\
b. Let $n>2$ and let $\mathcal{M}_{n-2}$ be the least inner model with $n-2$ Woodin cardinals, then $\mathbb Q$ adds a new $\Pi^1_n$ singleton over $\mathcal{M}_{n-2}$.
\\
c. Clause (b) is optimal in the following sense: If $\mathbb P \in \mathcal{M}_{n-2}$ is a Borel ccc forcing,  then $\mathbb P$ doesn't add a new $\Sigma^1_n$ singleton over $\mathcal{M}_{n-2}$.
\\
\\
\textbf{Proof}: Clauses (a) and (b) follow from Claim 24, with clause (a) using the fact that $\mathbb R \cap L$ is $\Sigma^1_2$ definable and clause (b) using Theorem 25. Clause (c) follows from Theorem 26. $\square$

\textbf{\large 9. Open Questions}{\large \par}

As our model doesn't sasitfy $AC_{\aleph_0}$,
it's natural to ask whether we can improve the result getting a model
of $AC_{\aleph_0}$ or even $DC$. Hopefully in {[}F1424{]} it will be shown that assuming
the existence of a measurable cardinal, we can get a model of $DC(\aleph_1)$.
This leads to the following question:

\textbf{Problem 1: }Can we improve the current result and get
a model of $DC$ without large cardinals?

A very recent work in preparation ([Sh1257]) answers this question in the affirmative.

As the current result gives measurability with respect to the ideal
$I_{\mathbf n,\aleph_1}$, it's natural to ask:

\textbf{Problem 2: }Can we get a similar result for the ideal
$I_{\mathbf n,\aleph_0}$?

This problem will be addressed in {[}GHS1097{]}.

\textbf{\large 10. References}{\large \par}

{[}Br{]} Joerg Brendle, Mad families and iteration theory, Logic and algebra, vol. 302 of Contemp. Math., pages 1-31. Amer. Math. Soc., Providence, RI, 2002.

{[}GHS1097{]} Mohammad Golshani, Haim Horowitz and Saharon Shelah. On the classification of definable
ccc forcing notions, arXiv:1610.07553

{[}F1424{]} Haim Horowitz and Saharon Shelah. Saccharinity with ccc:
Getting $DC(\omega_1)$ from a measurable cardinal, in preparation

{[}Je{]} Ronald Jensen, Definable sets of minimal degree, Studies in Logic and the Foundations of Mathematics, Volume 59, 1970, 122-128

{[}JuSh292{]} Haim Judah and Saharon Shelah. Souslin forcing. The
Journal of Symbolic Logic, 53: 1188-1207, 1988

{[}KeSo{]} Alexander S. Kechris and Slawomir Solecki. Approximation
of analytic by Borel sets and definable countable chain conditions.
Israel Journal of Mathematics, 89:343-356, 1995.

{[}KrSh859{]} Jakob Kellner and Saharon Shelah. Saccharinity. The
Journal of Symbolic Logic, 74: 1153-1183, 2011 

{[}KST{]} Jakob Kellner, Saharon Shelah and Anda R. Tanasie, Another ordering of the ten cardinal characteristics in Cichon's diagram, Comment. Math. Univ. Carolin. 60 (2019), no. 1, 61-95.

{[}Me1{]}  Diego A. Mejia, Anatomy of $\tilde{\mathbb E}$, arXiv:2402.04706

{[}Me2{]}  Diego A. Mejia, Template iterations with non-definable ccc forcing notions. Ann. Pure Appl. Logic, 166(11):1071-1109, 2015.

{[}MSW{]} Sandra Mueller, Ralf Schindler and W. Hugh Woodin, Mice from optimal deter- minacy hypotheses, Journal of Mathematical Logic. Volume 20, Issue Supp01, October 2020. 1950013.

{[}RoSh470{]} Andrzej Roslanowski and Saharon Shelah. Norms on possibilities
I: forcing with trees and creatures. Memoirs of the American Mathematical
Society, 141(671), 1999.

{[}RoSh628{]} Andrzej Roslanowski and Saharon Shelah. Norms on possibilities
II: More ccc ideals on $2^{\omega}$. Journal of Applied Analysis
3 (1997) 103-127.

{[}RoSh672{]} Andrzej Roslanowski and Saharon Shelah, Sweet \& sour
and other flavours of ccc forcing notions. Arch. Math. Logic 43 (2004)
583-663.

{[}Sh176{]} Saharon Shelah. Can you take Solovay's inaccessible away?
Israel Journal of Mathematics, 48:1-47, 1984.

{[}Sh1257{]} Saharon Shelah, Homogeneous forcing, in preparation

{[}Sh700{]} Saharon Shelah. Two cardinal invariants of the continuum
$\mathfrak d<\mathfrak a$ and FS linearly ordered iterated forcing.
Acta Math 192 (2004) 187-223 

{[}So1{]} Slawomir Solecki. Analytic ideals and their applications.
Annals of Pure and Applied Logic, 99:51-72, 1999.

{[}So2{]} Robert Solovay. A model of set theory in which every set
of reals is Lebesgue measurable. Annals of Mathematics. Second series
92 (1): 1-56

{[}St{]} John Steel, Projectively well-ordered inner models, Annals of Pure and Applied Logic 74 (1995) 77-104

{[}Za{]} Jindrich Zapletal, Forcing idealized. Cambridge Tracts in
Mathematics, vol. 174, Cambridge University Press, Cambridge, 2008.
MR MR2391923 (2009b:03002)

$\\$

(Haim Horowitz) Eintstein Institute of Mathematics

Edmond J. Safra campus, 

The Hebrew university of Jerusalem.

Givat Ram, Jerusalem, 91904, Israel.

E-mail address: haim.horowitz@mail.huji.ac.il

$\\$

(Saharon Shelah) Eintstein Institute of Mathematics

Edmond J. Safra campus, 

The Hebrew university of Jerusalem.

Givat Ram, Jerusalem, 91904, Israel.

Department of mathematics

Hill center - Busch campus, 

Rutgers, The state university of New Jersey.

110 Frelinghuysen road, Piscataway, NJ 08854-8019 USA

E-mail address: shelah@math.huji.ac.il

$\\$
\end{document}